\documentclass[aps, pre, twocolumn, superscriptaddress]{revtex4} 

\usepackage[switch]{lineno}
\usepackage{amsmath}
\usepackage{amssymb}
\usepackage{empheq}
\usepackage[parfill]{parskip}
\usepackage[active]{srcltx}
\usepackage{color}
\usepackage{array}
\usepackage{booktabs}
\usepackage{amsfonts}
\usepackage{dsfont}
\usepackage{graphicx}
\usepackage{hyperref}
\usepackage[dvipsnames]{xcolor}
\hypersetup{colorlinks,linkcolor={blue},citecolor={blue},urlcolor={blue}}  
\usepackage[english]{babel}
\usepackage[autostyle]{csquotes}
\usepackage{subcaption}
\usepackage{mwe}

\begin{document}
\setcitestyle{square}
\setlength{\parindent}{0.5cm}

\title{Attractive and repulsive interactions in the one-dimensional swarmalator model}

\author{Baoli Hao}
\affiliation{Department of Applied Mathematics, Illinois Institute of Technology, Chicago, IL 60616} 
\email{bhao2@hawk.iit.edu, mzhong3@iit.edu}

\author{Ming Zhong}
\affiliation{Department of Applied Mathematics, Illinois Institute of Technology, Chicago, IL 60616} 
\email{mzhong3@iit.edu}

\author{Kevin O'Keeffe}
\affiliation{Senseable City Lab, Massachusetts Institute of Technology, Cambridge, MA 02139} 
\email{Corresponding to: kevin.p.okeeffe@gmail.com}

\begin{abstract}

We study a population of swarmalators, mobile variants of phase oscillators, which run on a ring and have both attractive and repulsive interactions. This 1D swarmalator model produces several of collective states: the standard sync and async states as well as a novel splay-like ``polarized" state  and several unsteady states such as active bands or swirling. The  model's simplicity allows us to describe some of the states analytically. The model can be considered as a toy model for real-world swarmalators such as vinegar eels and sperm which swarm in quasi-1D geometries.
\\
\noindent \\
\end{abstract}

\maketitle

\section{Introduction}
Large groups of coupled oscillators have been used to model diverse phenomena \cite{winfree2001geometry,kuramoto2003chemical,strogatz2000kuramoto}. Early studies considered oscillators with no spatial embedding; they were coupled all-to-all with uniform strength \cite{kuramoto2003chemical}. Later this simplification was relaxed and oscillators were given different arrangements in space. Rings of oscillators were studied, which produced splay states and chimeras \cite{ermentrout1985behavior,sakaguchi1987local,wiley2006size,abrams2004chimera}, and lattices, which produced vortices and spirals \cite{strogatz1988collective,ottino2016frequency,lee2010vortices}.

New work \cite{o2017oscillators,o2022collective,yoon2022sync} considers the next step in this sequences of generalizations: it considers oscillators which are free to move around in space — oscillators which sync and swarm. Swarmalators, short for swarming oscillators, couple their internal and external degrees of freedom bidirectionally: Their movements depend on their phases, just as their phases depend on their movements. With a view to explaining the behavior of biological microswimmers \cite{yang2008cooperation,riedel2005self,quillen2021metachronal,quillen2022fluid,tamm1975role,verberck2022wavy,belovs2017synchronized,peshkov2022synchronized}, chemical micromotors \cite{yan2012linking,hwang2020cooperative,zhang2020reconfigurable,bricard2015emergent,zhang2021persistence,manna2021chemical,li2018spatiotemporal,chaudhary2014reconfigurable}, and other system which both sync and swarm \cite{hrabec2018velocity,tan2021development,petrungaro2019information,gengel2023physics,ceron2023programmable} researchers have studied swarmalators with pinning \cite{sar2023pinning,sar2023solvable,sar2023swarmalators}, local coupling \cite{lee2021collective}, stochastic coupling \cite{schilcherswarmalators}, delayed coupling \cite{blum2022swarmalators}, external forcing \cite{lizarraga2020synchronization}, phase frustration \cite{lizarraga2023synchronization} and other effects \cite{o2018ring,jimenez2020oscillatory,ha2021mean,degond2022topological,sar2022dynamics,degond2023topological,ceron2023diverse,hong2023swarmalators,sync_saka,gerew2023concurrent,ha2023fast,yadav2023exotic,smith2023swarmalators}. Applications of swarmalators to robotics  have also been considered \cite{barcis2019robots,barcis2020sandsbots,rinnermultidrone,monaco2019cognitive,monaco2020cognitive,ceron2023towards}. 

Here we add to this young literature by studying swarmalators with a mix of attractive and repulsive interactions. Mixed sign couplings (mixed sign coupling meaning a mixture of positive/attractive and negative/repulsive couplings) are common in systems of regular oscillators, for example in neurons which have both attractive and inhibitory couplings \cite{borgers2003synchronization}. We suspect they are also common in systems of swarmalators. Active colloids, for instance, have hydrodynamic interactions which can switch from being attractive to repulsive as the relative orientation between particles changes \cite{liebchen2021interactions}. A theoretical understandings of how mixed sign couplings changes the phenomenology of swarmalators in lacking. A first step towards filling in this gap was recently taken \cite{hong2021coupling} by studying swarmalators which move around in 2D. New states were found, but were unfortunately analytically intractable. Hence, we restrict the swarmalators movements to a one-dimensional (1D) ring with a view to making the analysis simpler (this 1D model may also be derived from the 2D swarmalator model, and in that sense captures the essence of 2D swarmalator phenomena). This 1D model could also be used for real-world modeling purposes, since sperm, Janus colloids, and other natural swarmalators are often confined to quasi-1D rings-like geometries \cite{bau2015worms,yuan2015hydrodynamic,ketzetzi2021activity,creppy2016symmetry,aihara2014spatio} and likely have a mix of attractive and repulsive interactions in some settings. Our main findings are a variety of new collective states which we describe with a mix of theory and numerics.
 
\section{Model}
The 1D swarmalator model \cite{o2022collective,yoon2022sync} we study is
\begin{align}
    \dot{x_i} &= \nu+\frac{J_i}{N} \sum_j^N \sin(x_j - x_i) \cos(\theta_j - \theta_i) \label{eom-x} \\
    \dot{\theta_i} &= \omega+\frac{K_i}{N}  \sum_j^N  \sin(\theta_j - \theta_i ) \cos(x_j - x_i ) \label{eom-theta}
\end{align}
\noindent
where $(x_i, \theta_i) \in (\mathbb{S}^1, \mathbb{S}^1)$ are the position and phase of the $i$-th swarmalator for $i = 1, \dots, N$, and 
$\nu$ and $\omega$ represent the constant natural frequency and $(J_i,K_i)$ are the associated couplings constants. Note $\mathbb{S}$ denotes the unit circle. For simplicity, we set $J_i = J = 1$ and draw the phase coupling $K$ from 
\begin{align}
h(K) &= p\delta(K - K_p) + q \delta(K - K_n)
\label{hKgJ}
\end{align}
where $p+q=1$. We see a fraction $p$ of the swarmalators have positive couplings $K_p > 0$, and the remaining $q = 1-p$ have negative coupling $K_n<0$. This mix of positive and negative coupling has been studied before in the regular Kuramoto model \cite{hong2011kuramoto}, where the oscillators with $K_p$ were called `conformists', since positive coupling tends to synchronize oscillators (in that sense the oscillators `conform`) and those with negative coupling $K_n$ tend to anti-synchronize (and in that sense are contrarian). We will use the same terminology here. For simplicity, we set $J = 1$ without loss of generality by rescaling time which leaves a model with three parameters $(p,K_p, K_n)$. The model with identical couplings $(J,K)$ \cite{o2022collective} and couplings of form $(J_j, K_j)$ \cite{o2022swarmalators} were previously studied; this $(J_i,K_i)$ coupling study is a natural generalization of these works. The 1D swarmalator model may also be derived from the 2D swarmalator model \cite{o2017oscillators} (see Appendix in \cite{o2022collective}). 

\section{Numerics}
Numerical experiments were performed by us to explore the behavior of our model. We used Matlab's ODE solver ``ode45" to run our simulations. The swarmalators are initially positioned in $[0,2\pi]$ and their initial phases were drawn in the same domain, both uniformly at random. We studied various parameters $(p,K_p,K_n)$ and observed seven collective states. Four of these are static, in the sense the individual swarmalators are ultimately stationary in both $x$ and $\theta$. In contrast, the remaining three were unsteady. 

We used three order parameters to catalog the states: the rainbow order parameters $W_{\pm}$ used in previous studies of swarmalators \cite{o2017oscillators} and the mean velocity $V$. Their definitions are
\begin{subequations}
    \begin{align}
         & W_{\pm} =  S_{\pm} e^{i \phi_{\pm}} := \frac{1}{N} \sum_{j=1}^N e^{i (x_j\pm \theta_j)} \quad (i=\sqrt{-1}), \label{eq:order parameter}\\
          &V := \frac{1}{N} \sum_{j=1}^N \langle| \dot{x}_j |\rangle_t. \label{eq:mean_velocity}
    \end{align}
\end{subequations}
The magnitudes $S_{\pm}$ where $0 \leq S_{\pm} \leq 1$ measure the amount of space-phase correlation. When $x_i$ and $\theta_i$ are uncorrelated, the order parameters take minimal values $S_{\pm} = 0$. When $x_i$ and $\theta_i$ are perfectly correlated, however, one is maximal, the other minimal $(S_+, S_-) = (1,0)$. The opposite happens when $x_i, \theta_i$ are anti-correlated $x_i = - \theta_i + C$, $S_+ = 0, S_- = 1$. The symmetry in our model means that perfect correlation and anti-correlation occur equally, so we instead define $S_{max}, S_{min} = \max{S_{\pm}}, \min{S_{\pm}}$ which eliminates this degeneracy (in the sense that $S_{max}$ is always $1$ and $S_{min} = 0$ for the both the correlated and anti-correlated cases). Finally, we note that when the positions are fully sync'd $x_i = C_1$ and the phase are fully sync'd $\theta_i = C_2$, then both $S_{\pm}$ are maximal simultaneously $S_+ = S_- = 1$. 

We next discuss each of the collective states. We recommend viewing Supplementary Movie 1 at this point, which shows all the states at once \cite{sm}. Having this visual in mind will be helpful when reading the verbal descriptions and associated figures. 
We also provide a github link \footnote{\href{https://github.com/Khev/swarmalators/tree/master/1D/on-ring/mixed-coupling/ki}{https://github.com/Khev/swarmalators/tree/master/1D/on-ring/mixed-coupling/ki}} to code to simulate the model, which can be helpful to get a deeper understanding of the states.

\begin{figure}
        \centering
        \begin{subfigure}[b]{0.225\textwidth}
           \centering
            \includegraphics[width=\textwidth]{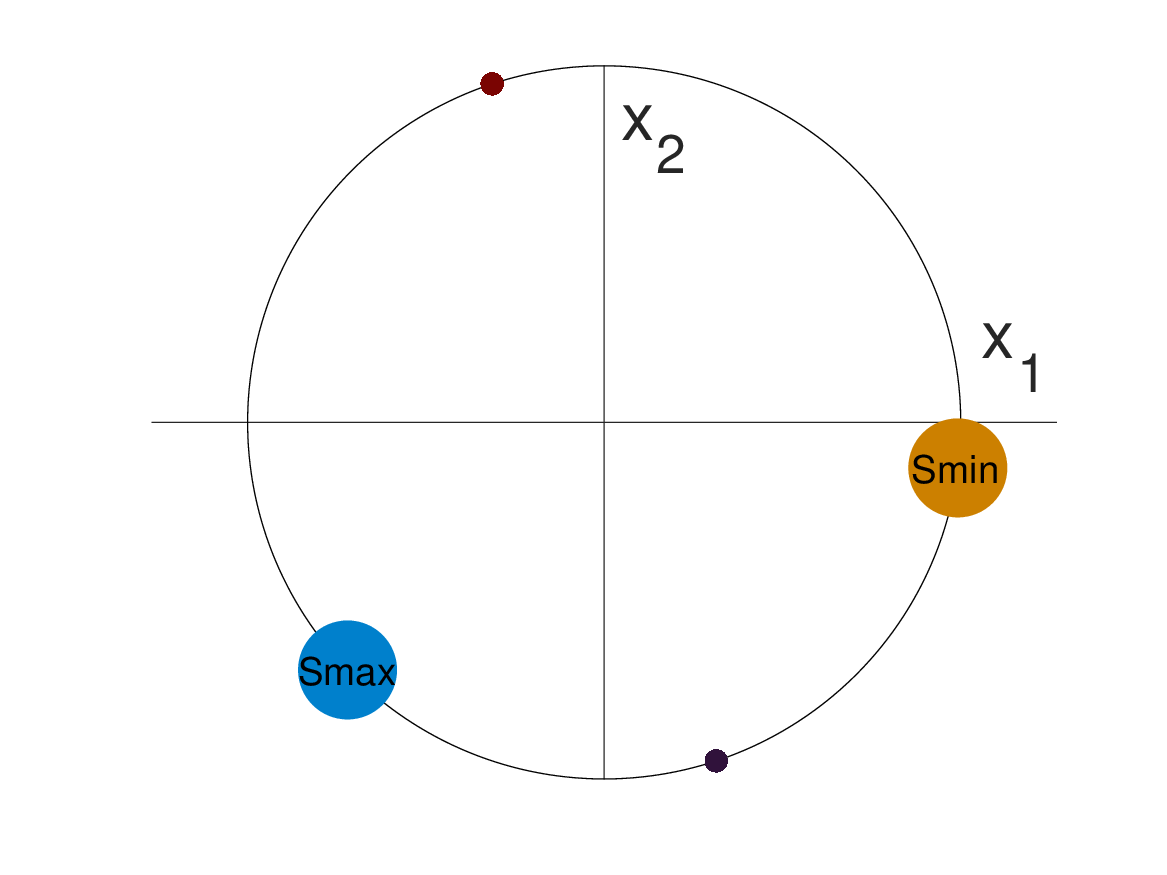}
            \caption[]%
            {{\small Static Synchrony}}    
            \label{fig:Static Synchrony}
        \end{subfigure}
        \hfill
        \begin{subfigure}[b]{0.225\textwidth}  
            \centering
            \includegraphics[width=\textwidth]{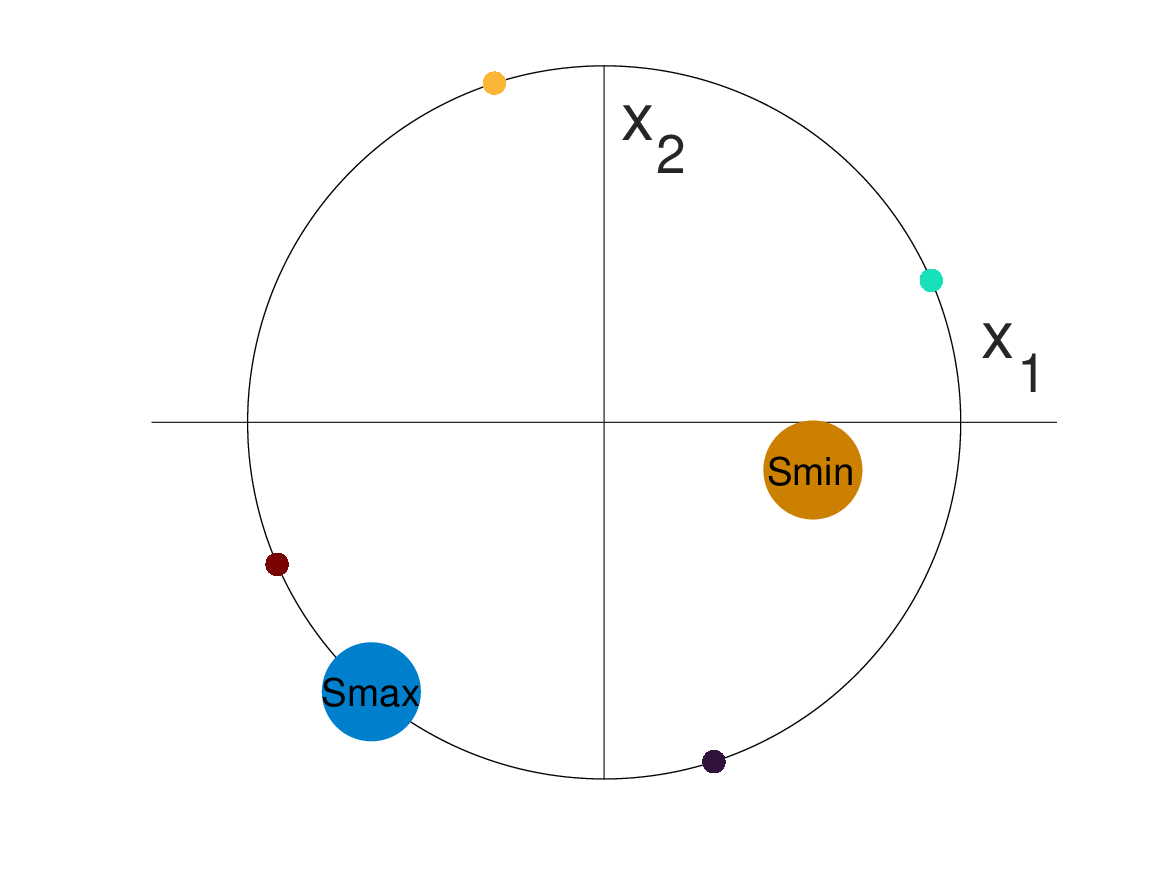}
            \caption[]%
            {{\small Polarized State}}    
            \label{fig:Polarized State}
        \end{subfigure}
      \vskip\baselineskip
        \begin{subfigure}[b]{0.225\textwidth}   
            \centering
            \includegraphics[width=\textwidth]{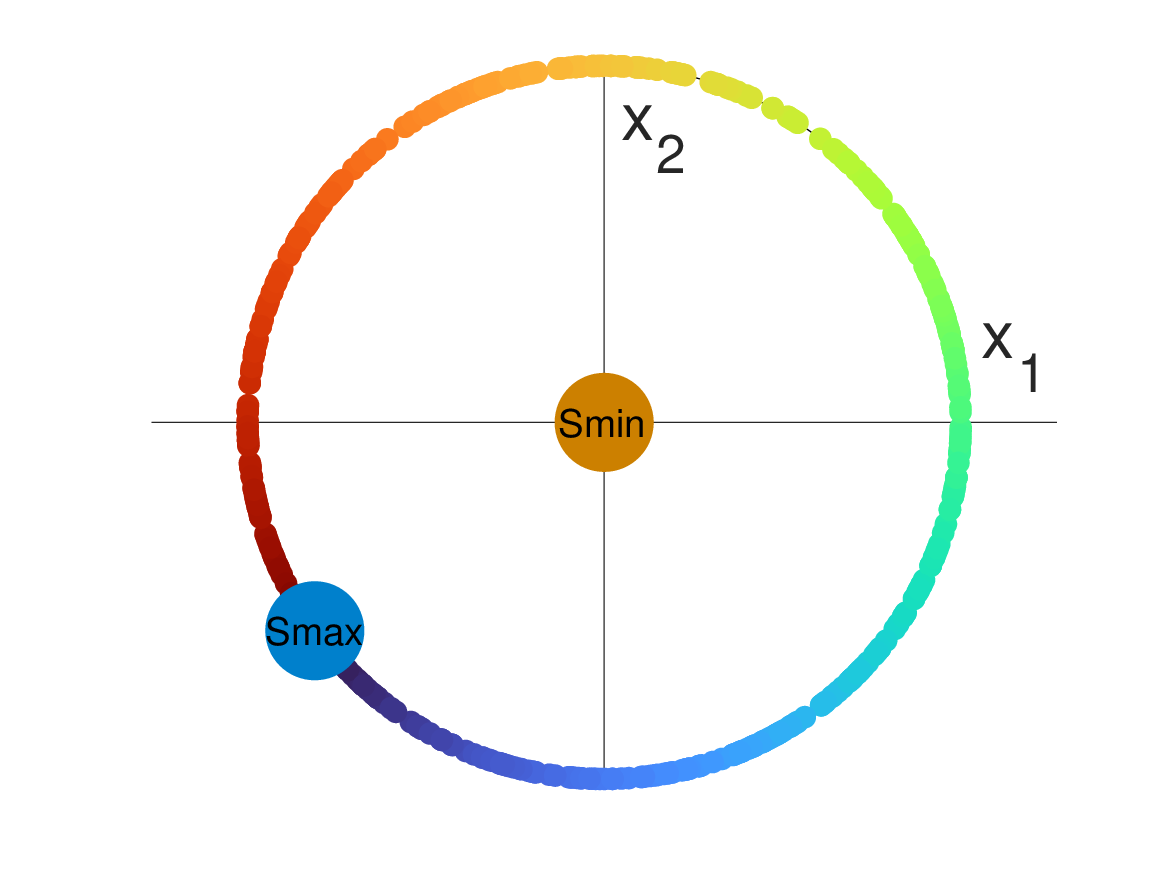}
            \caption[]%
            {{\small Static Phase Wave}}    
            \label{fig:Static Phase Wave}
        \end{subfigure}
        \hfill
        \begin{subfigure}[b]{0.225\textwidth}   
           \centering
            \includegraphics[width=\textwidth]{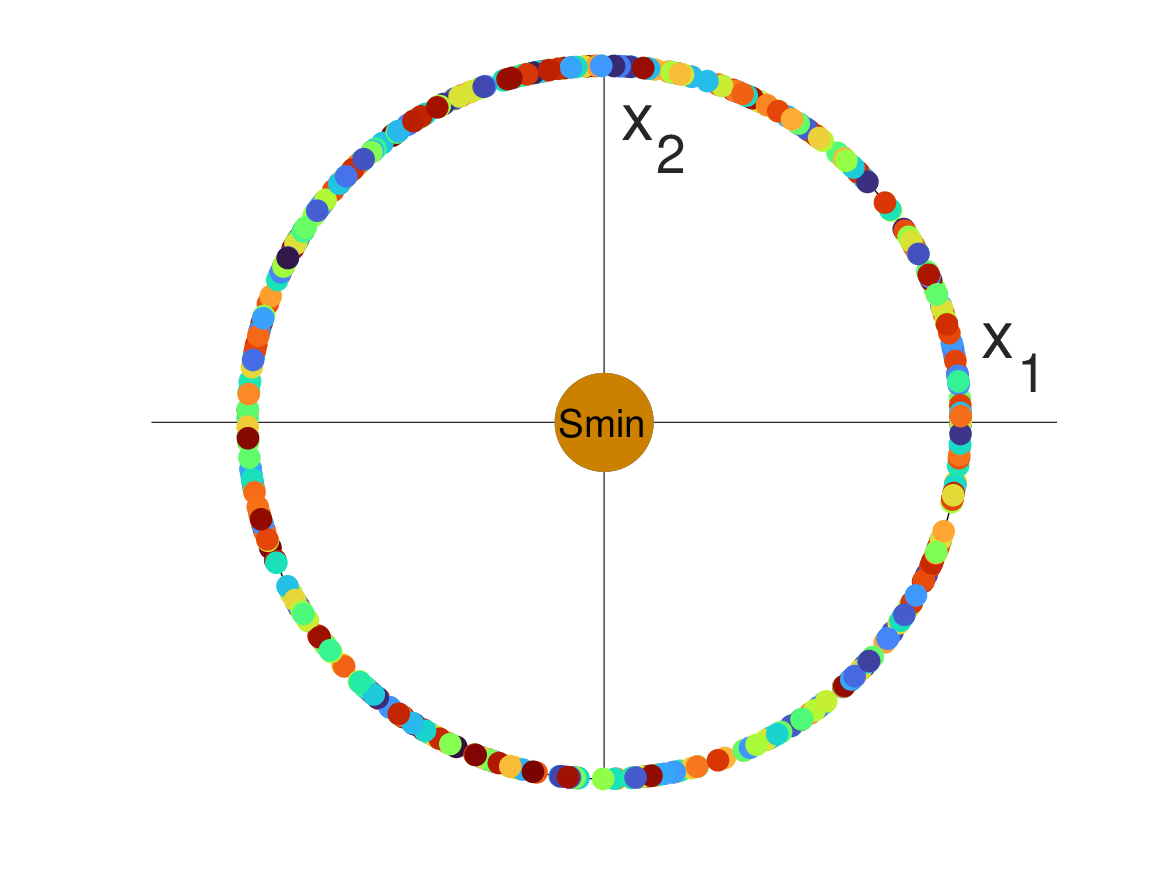}
            \caption[]%
            {{\small Static Asynchrony}}    
            \label{fig:Static Asynchrony}
        \end{subfigure}
        \caption[ Stationary States ]
        {\small \textbf{Stationary Collective States.} Scatter plots of four stationary states in the $(x_1, x_2)$ plane, where $(x_1,x_2)=(\cos{x},\sin{x})$ and the swarmalators are colored in terms of their phases. Simulations were run with $N=500$ swarmalators  for variable numbers of time units $T$ and step size $dt=0.1$.  (a) Static sync state for $(J,K_n,K_p,p)=$ $(1,-0.5,0.5,1)$ and $T=100$.
        (b) Polarized state for $(J,K_n,K_p,p)=(1,-0.5,0.5,0.8)$ and $T=1000$. (c) Static phase wave state for $(J,K_n,K_p,p)=$ $(1,-0.5,0.5,0.2)$ and $T=100$. (d) Static async state for $(J,K_n,K_p,p)=$ $(1,-3,0.5,0.1)$ and $T=100$.} 
        \label{Scatter Plots}
    \end{figure}

\begin{figure}
        \centering
        \begin{subfigure}[b]{0.225\textwidth}
           \centering
            \includegraphics[width=\textwidth]{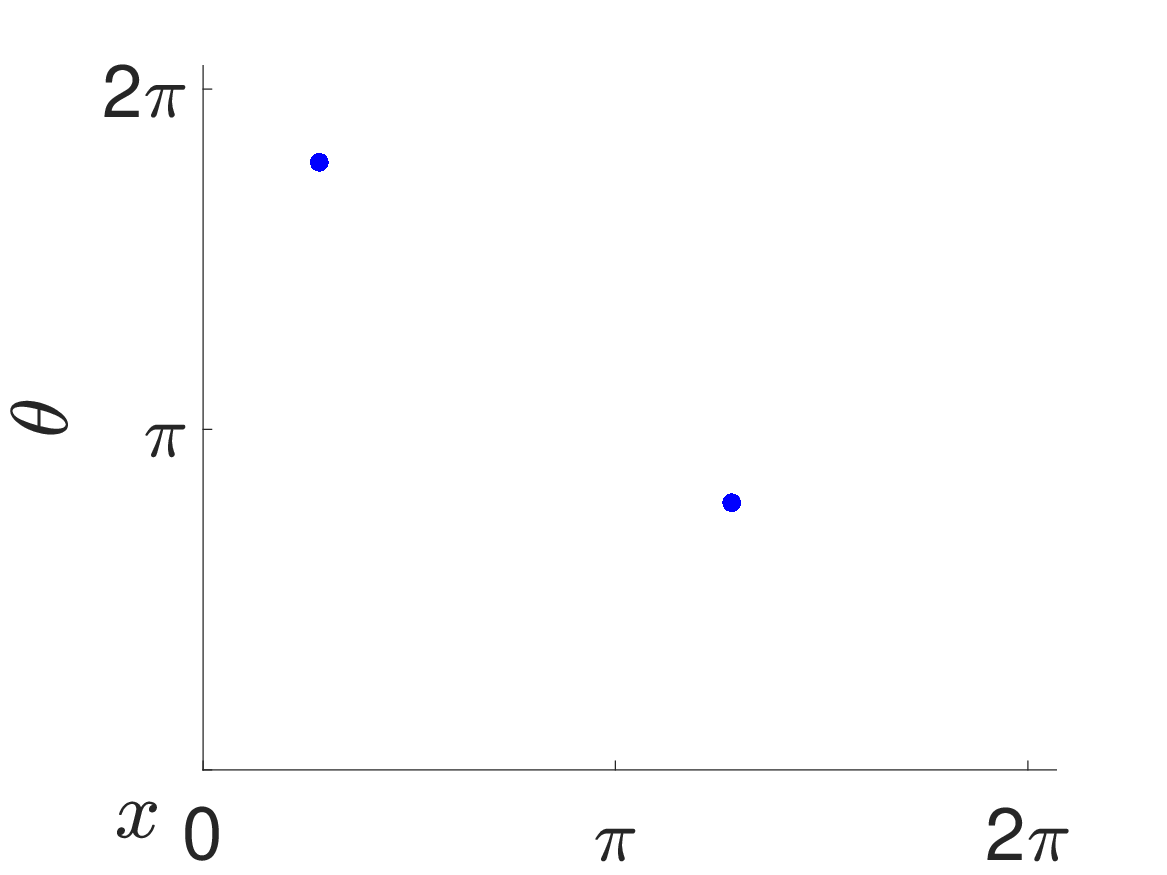}
            \caption[Static Synchrony]%
            {{\small Static Synchrony}}    
            \label{fig:xtheta_Static Synchrony}
        \end{subfigure}
        \hfill
        \begin{subfigure}[b]{0.225\textwidth}  
            \centering
            \includegraphics[width=\textwidth]{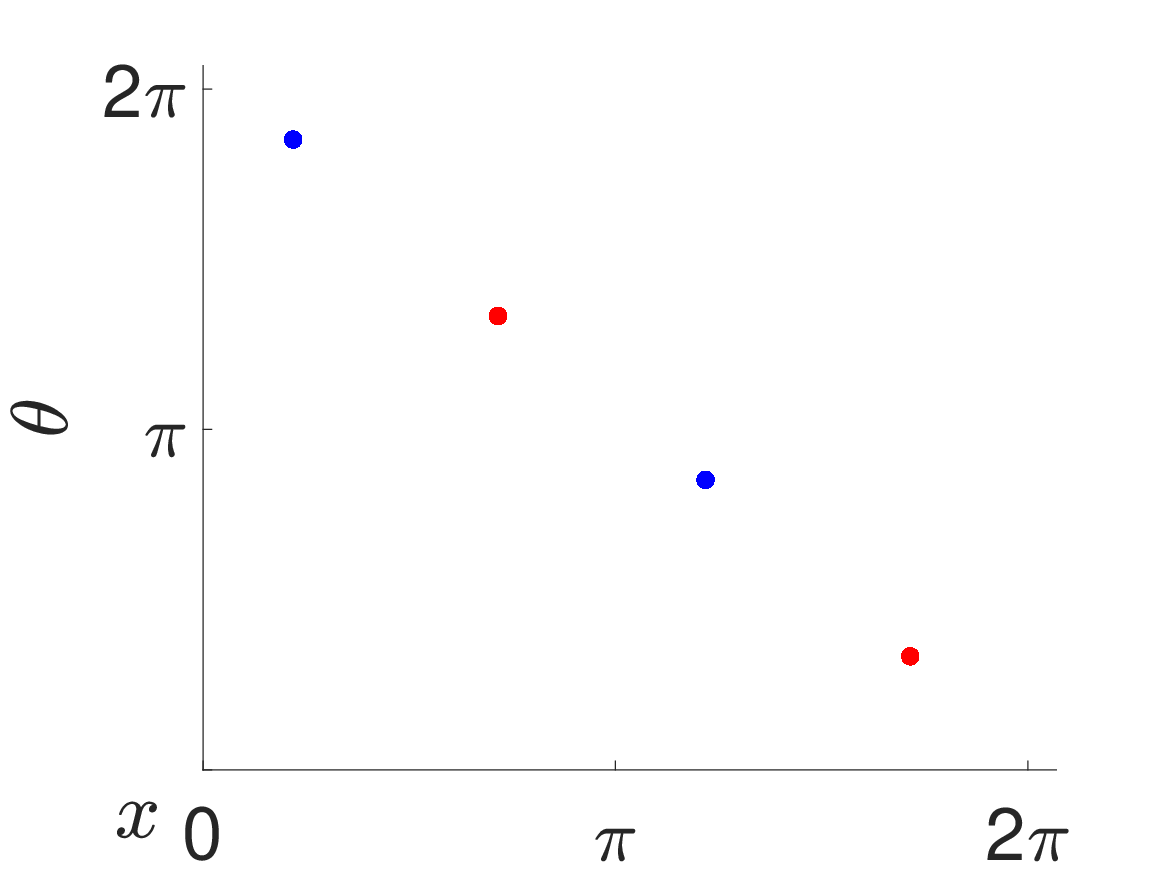}
            \caption[Polarized State]%
            {{\small Polarized State}}    
            \label{fig:xtheta_Polarized State}
        \end{subfigure}
      \vskip\baselineskip
        \begin{subfigure}[b]{0.225\textwidth}   
            \centering
            \includegraphics[width=\textwidth]{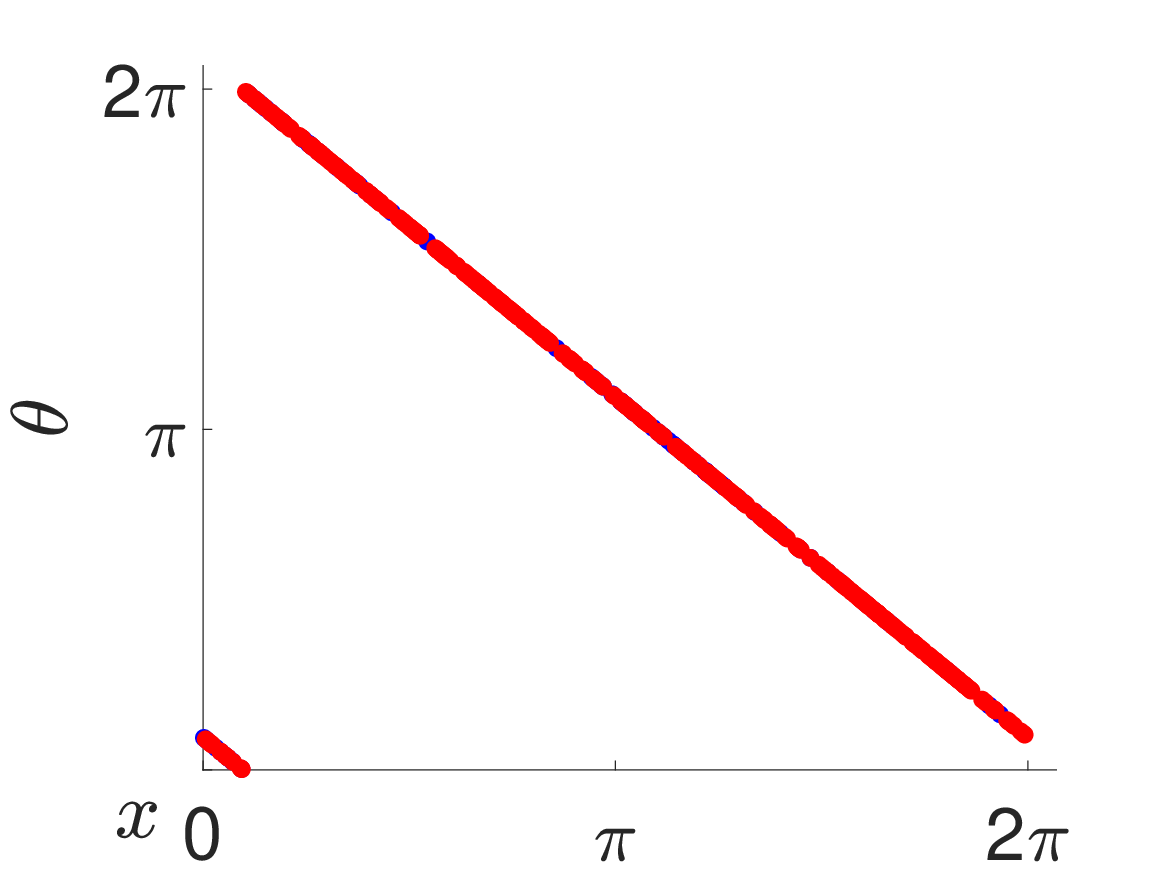}
            \caption[Static Phase Wave]%
            {{\small Static Phase Wave}}    
            \label{fig:xtheta_Static Phase Wave}
        \end{subfigure}
        \hfill
        \begin{subfigure}[b]{0.225\textwidth}   
           \centering
            \includegraphics[width=\textwidth]{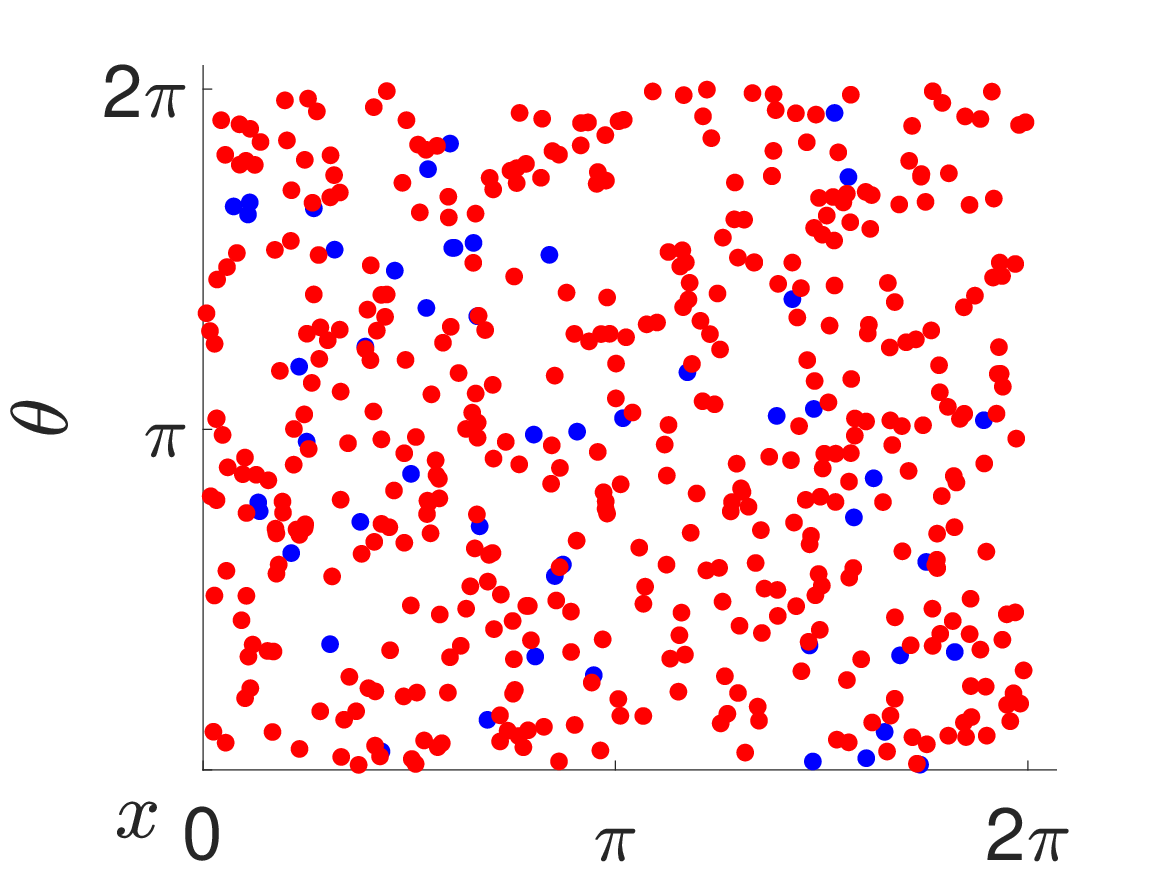}
            \caption[]%
            {{\small Static Asynchrony}}    
            \label{fig:xtheta_Static Asynchrony}
        \end{subfigure}
        \caption[ Scatter Plots ]
        {\small \textbf{Scatter plots in $(x,\theta)$ space.} Distributions in $(x, \theta)$ space corresponding to different states. Simulations were run with $N=500$ swarmalators  for variable numbers of time units $T$ and step size $dt=0.1$. Swarmalators coupling with $K_p$ and $K_n$ are presented as blue dots and red dots respectively. (a) Static sync state for $(J,K_n,K_p,p)=$ $(1,-0.5,0.5,1)$ and $T=100$.
        (b) Polarized state for $(J,K_n,K_p,p)=(1,-0.5,0.5,0.8)$ and $T=500$. (c) Static phase wave state for $(J,K_n,K_p,p)=$ $(1,-0.5,0.5,0.2)$ and $T=100$. (d) Static async state for $(J,K_n,K_p,p)=$ $(1,-3,0.5,0.1)$ and $T=100$.} 
        \label{fig:Scatter Plots}
    \end{figure}
    
\begin{figure}
    \centering
    \includegraphics[width= 1.05 \columnwidth]{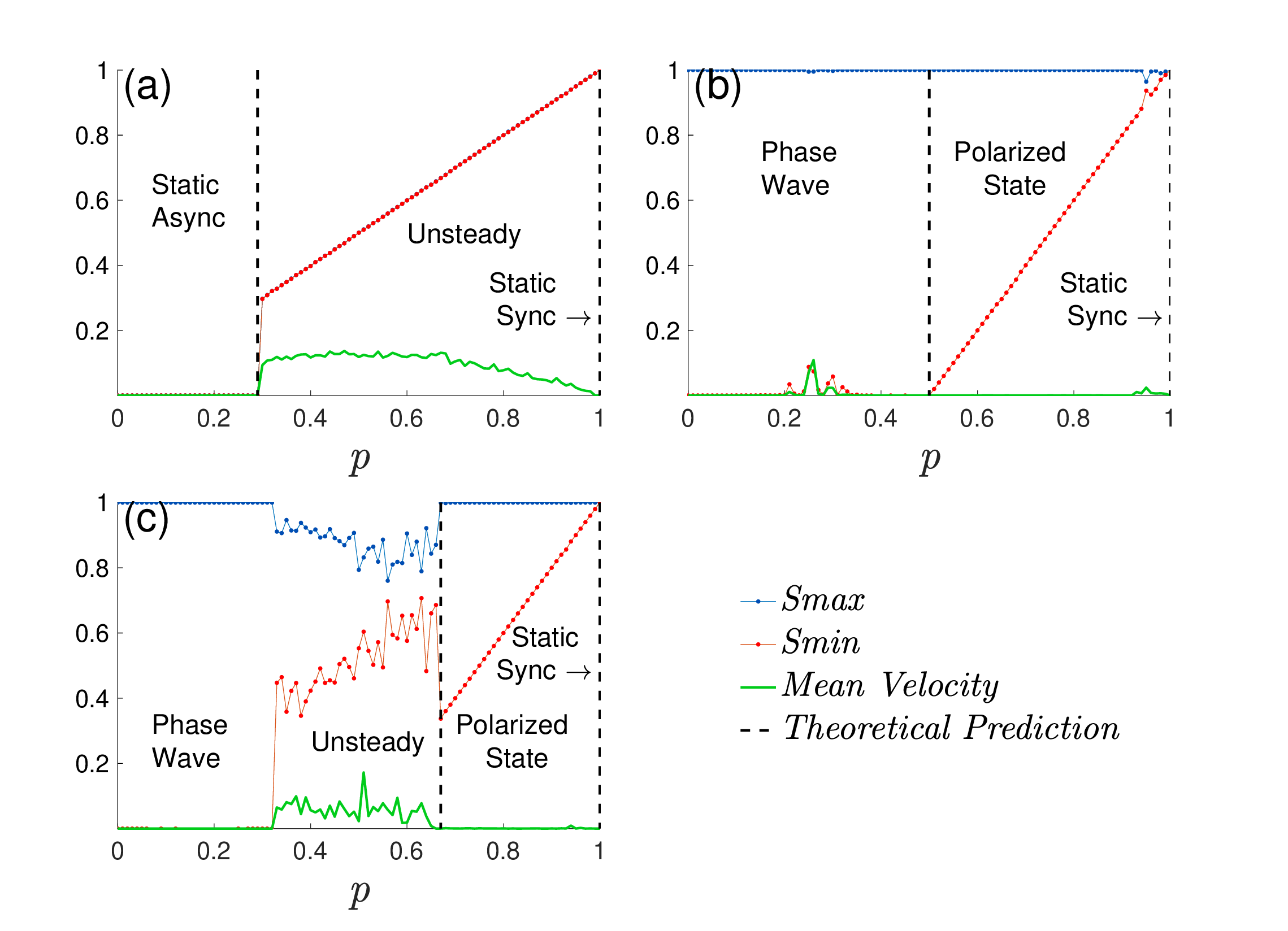}  
    \caption[ order parameters plot ]
    {\small \textbf{Order parameters and averaged velocity for different coupling distributions.} Asymptotic behavior of the order parameters $S_{max}$ := $\max(S_{+},S_{-}$) (blue dots) and $S_{min}$ := $\min(S_{+},S_{-}$)(red dots) versus $p$ for other parameters $(J,K_p,N,T,dt)$ = $(1,0.5,500,1000,0.1)$. (a) It shows the transition from static async to unsteady state with $p$ varying from 0 to 1 when $K_n = -2$. (b) It shows the transition from phase wave to polarized state with $p$ varying from 0 to 1 when $K_n=-0.8$. (c) It shows the transitions from phase wave to unsteady state and then polarized state when $K_n=-0.25$.  Each data point represents the average of last $10\%$ realizations.} 
        \label{fig:order parameters}
\end{figure}
   \begin{figure}
            \centering
            \includegraphics[width=0.5\textwidth]{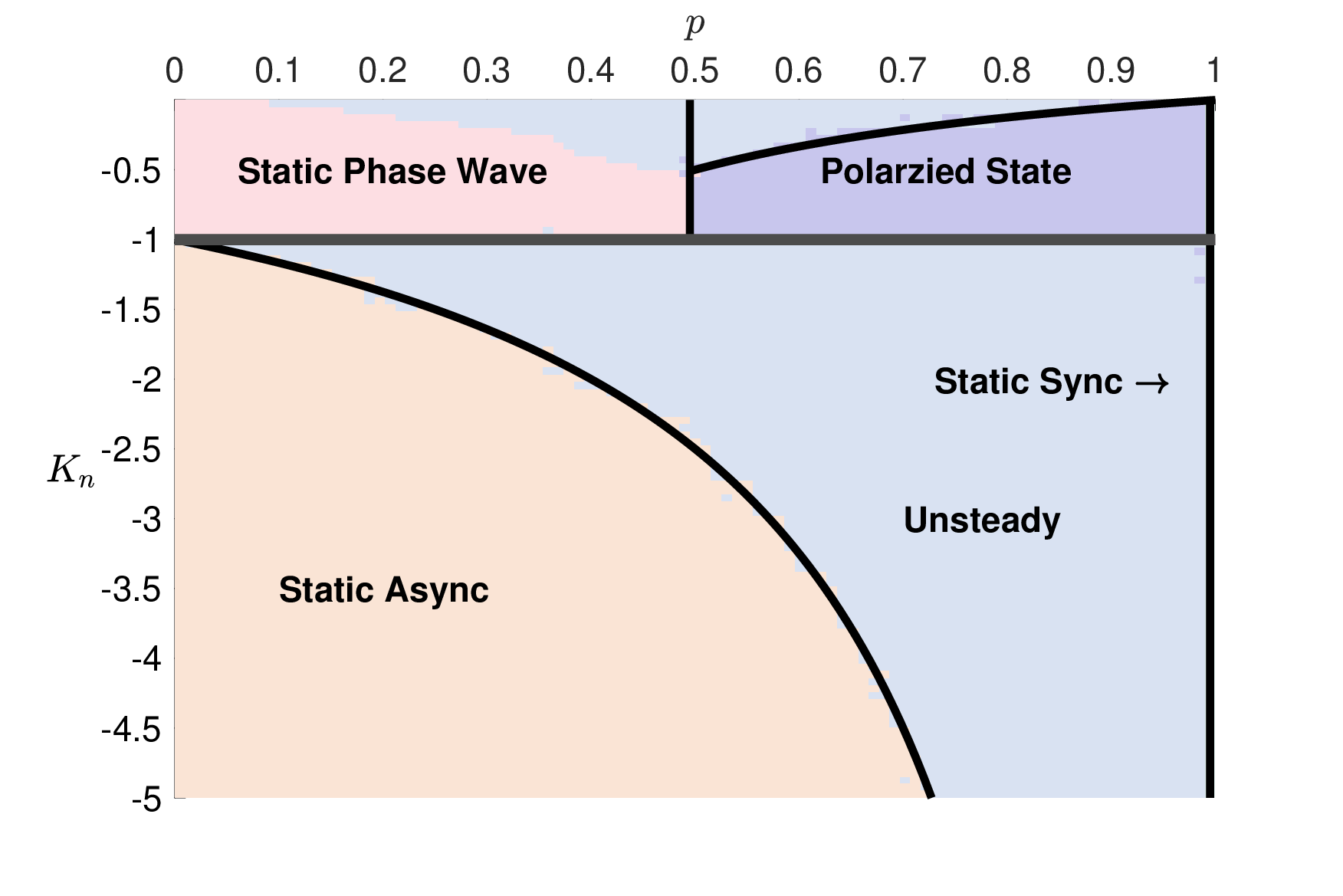}
            \caption[ phase Diagram ]
    {\small \textbf{Phase Diagram in $(P,K_n)$ Plane with fixed $K_p$=0.5.} Each state is indicated by a distinct color. The black curves and lines represent the theoretical predictions. Parameters in simulation we used are $(J,N,T,dt)=(1,5000,1000,0.1)$.} 
            \label{fig:meshgrid}
        \end{figure}

   \begin{figure*}
   
        \centering
      
         \begin{subfigure}[b]{0.25\textwidth}  
            \centering
            \includegraphics[width=\textwidth]{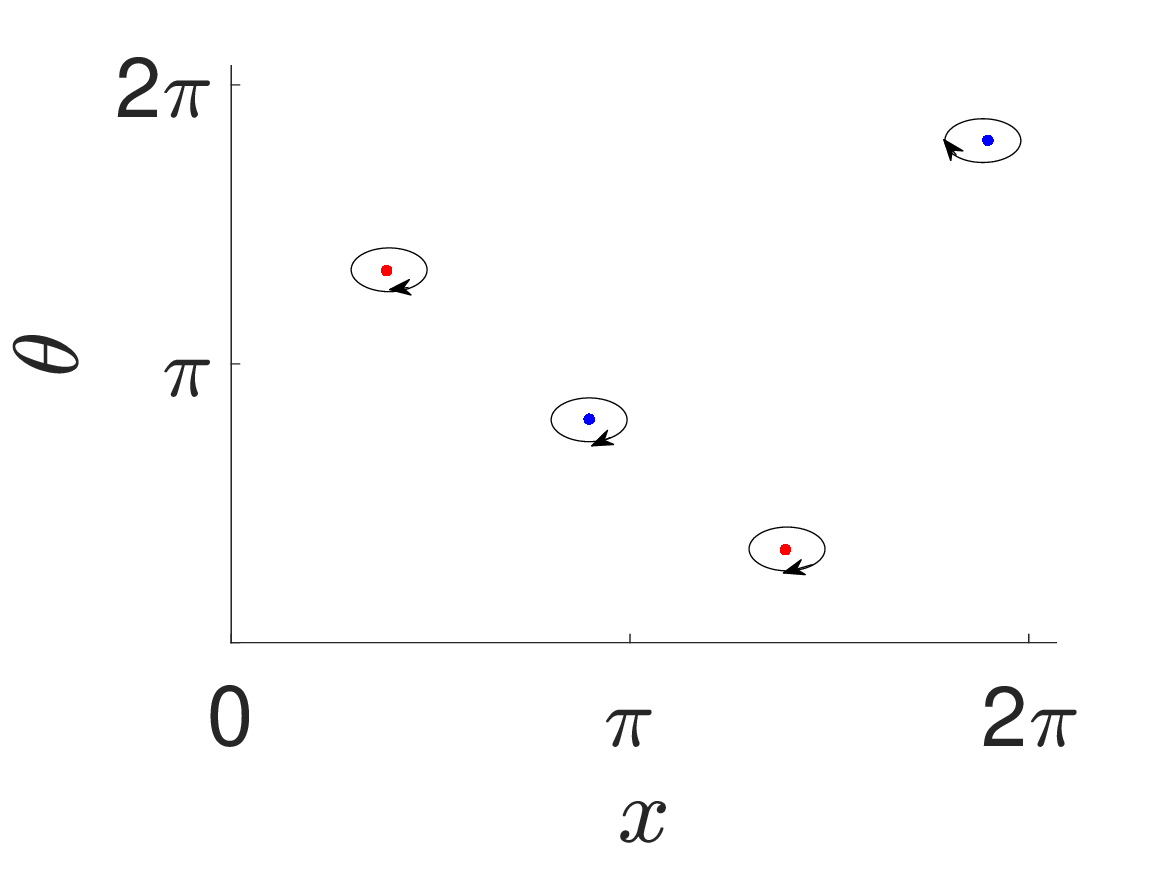}
            \caption[]%
            {{\small Breathing Polarized State}}    
            \label{fig:breathps}
        \end{subfigure}
        \hfill
        \begin{subfigure}[b]{0.25\textwidth}  
            \centering
            \includegraphics[width=\textwidth]{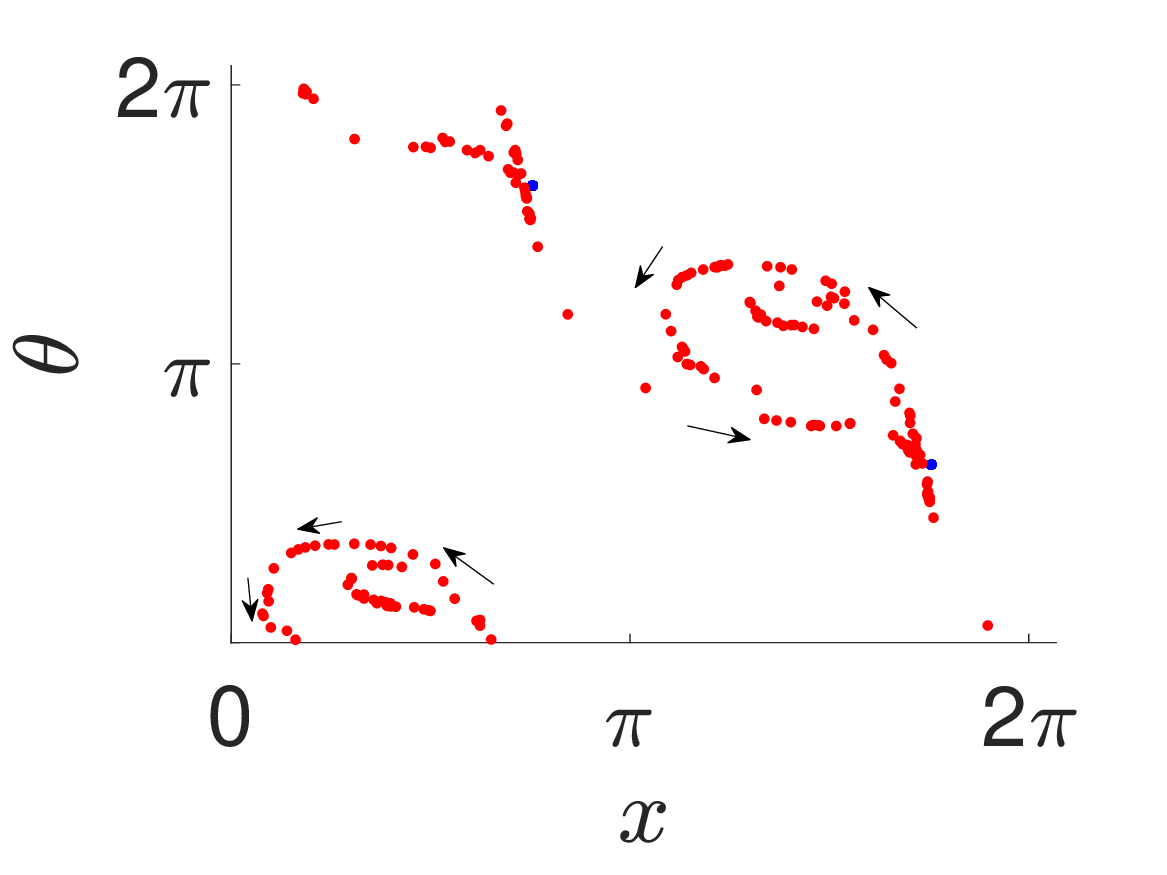}
            \caption[]%
            {{\small Swirling}}    
            \label{fig:swirling}
        \end{subfigure}
        \hfill
          \begin{subfigure}[b]{0.25\textwidth}
           \centering
            \includegraphics[width=\textwidth]{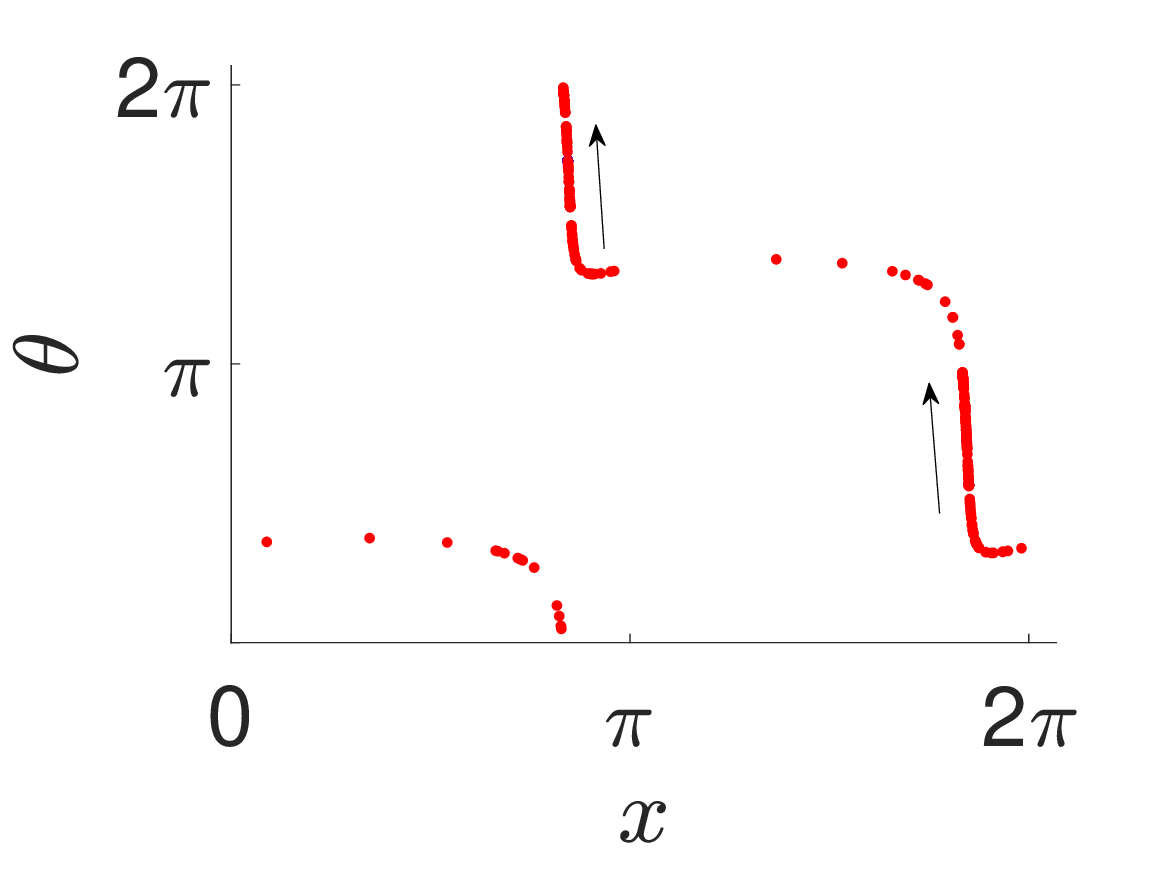}
            \caption[]%
            {{\small Active Bands }}    
            \label{fig:active_bands}
        \end{subfigure}
      \vskip\baselineskip
        \begin{subfigure}[b]{0.25\textwidth}  
            \centering
            \includegraphics[width=\textwidth]{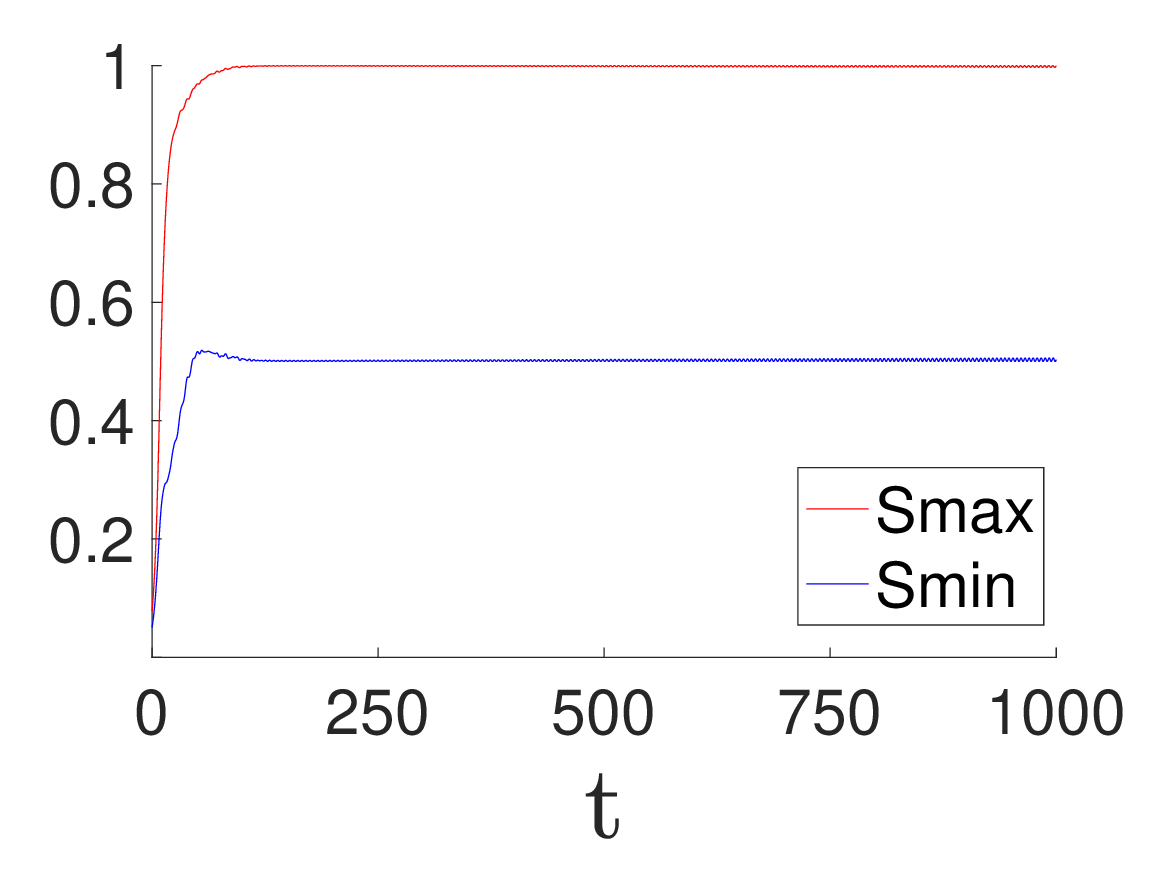}
            \caption[]%
            {{\small }}    
            \label{fig:xy_breathps}
        \end{subfigure}
        \hfill
        \begin{subfigure}[b]{0.25\textwidth}   
           \centering
            \includegraphics[width=\textwidth]{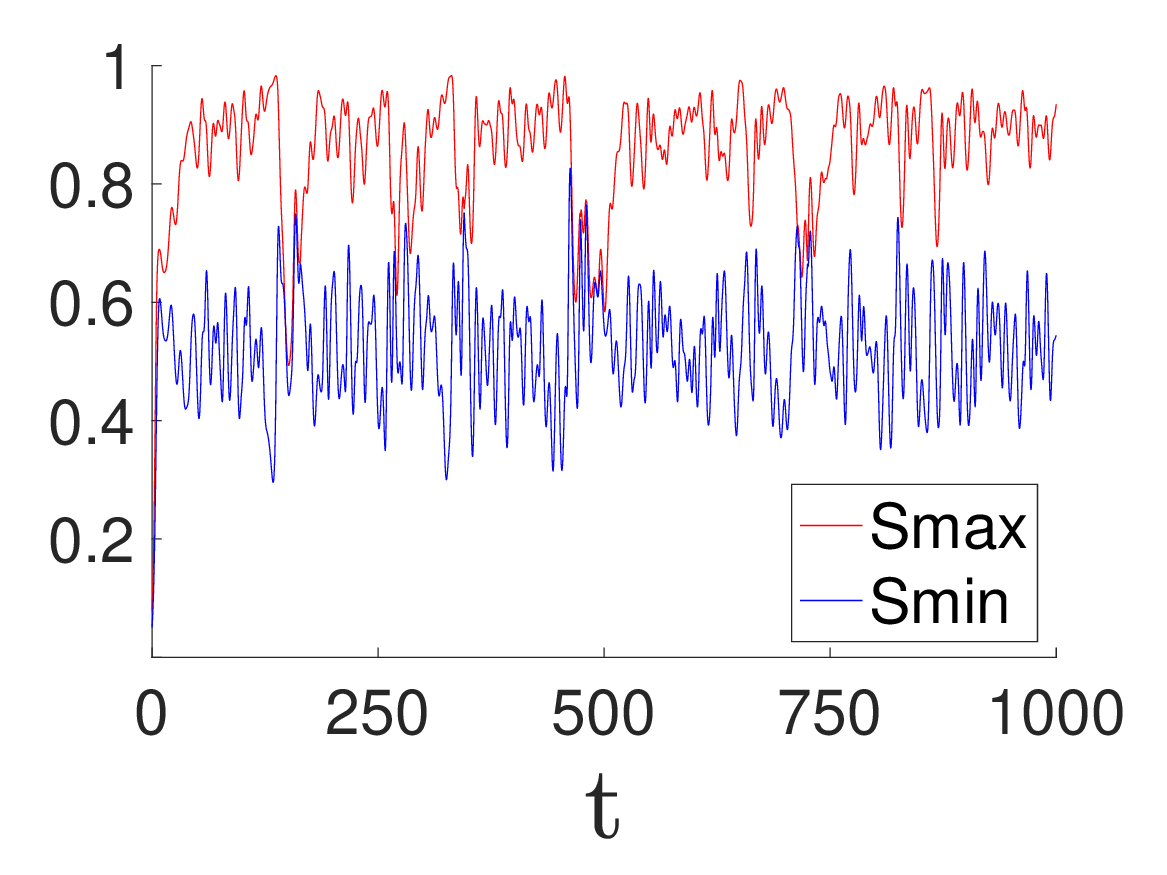}
            \caption[]%
            {{\small }}    
            \label{fig:xy_swirling}
        \end{subfigure}
        \hfill
         \begin{subfigure}[b]{0.25\textwidth}   
            \centering
            \includegraphics[width=\textwidth]{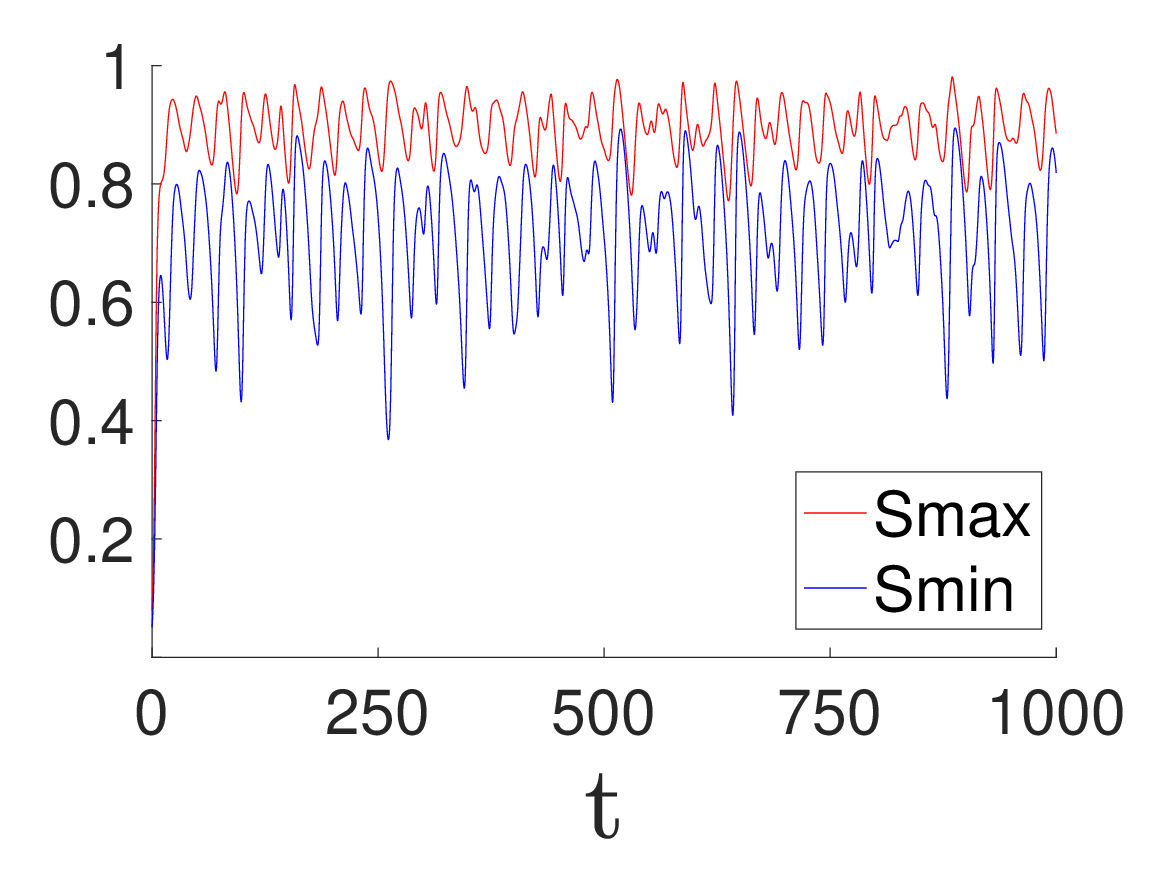}
            \caption[]%
            {{\small }}    
            \label{fig:xy_active_band}
        \end{subfigure}
        \vskip\baselineskip
        \begin{subfigure}[b]{0.25\textwidth}  
            \centering
            \includegraphics[width=\textwidth]{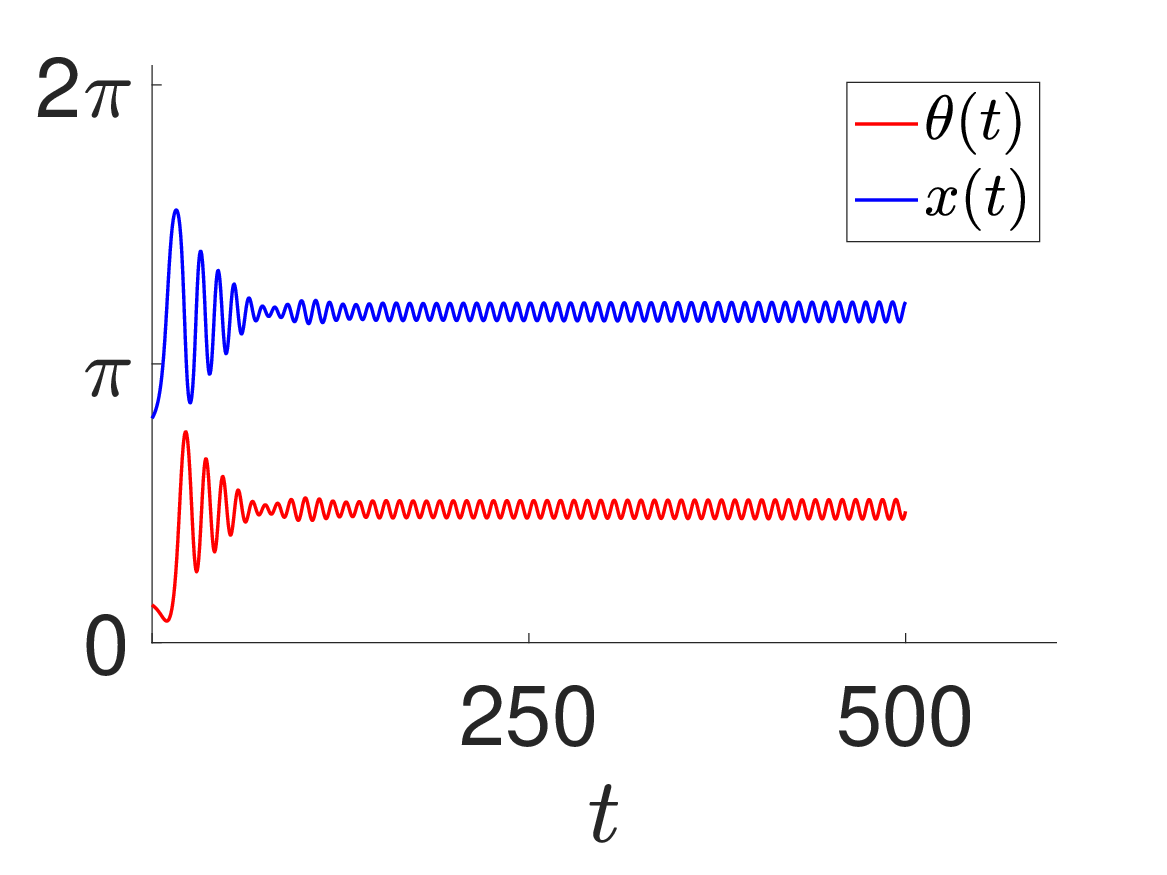}
            \caption[]%
            {{\small }}    
            \label{fig:1_breathps}
        \end{subfigure}
        \hfill
        \begin{subfigure}[b]{0.25\textwidth}   
           \centering
            \includegraphics[width=\textwidth]{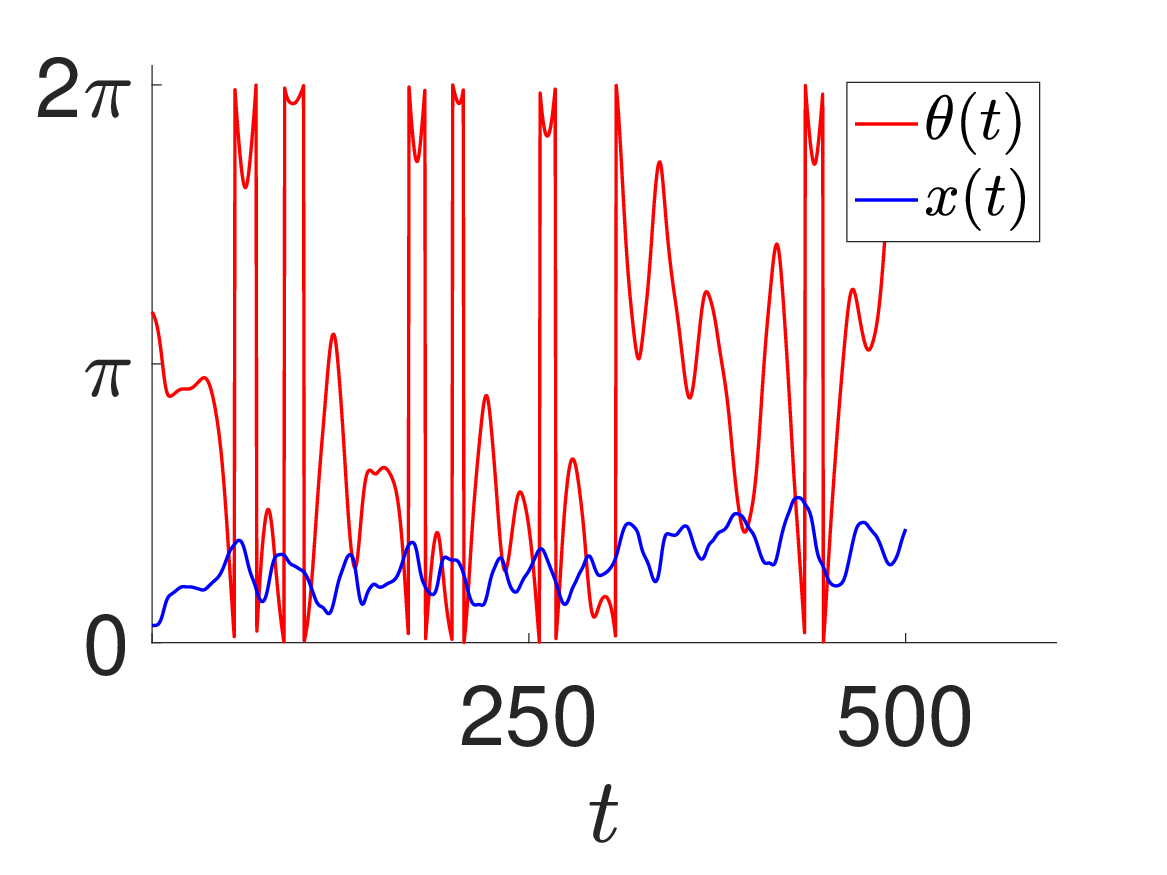}
            \caption[]%
            {{\small }}    
            \label{fig:1_swirling}
        \end{subfigure}
        \hfill
         \begin{subfigure}[b]{0.25\textwidth}   
            \centering
            \includegraphics[width=\textwidth]{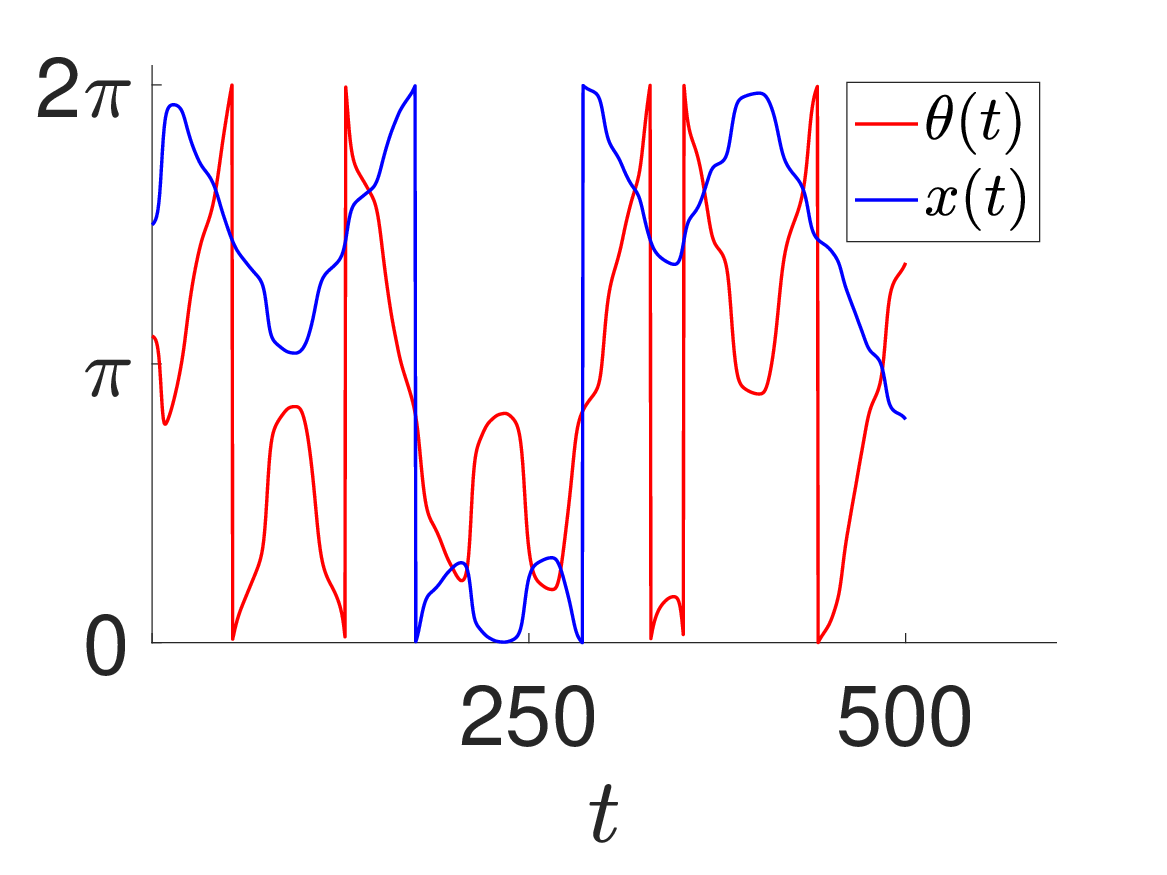}
            \caption[]%
            {{\small }}    
            \label{fig:1_active_band}
        \end{subfigure}
        
        \caption[ ]
        {\small \textbf{Unsteady Collective States.}} Top row: scatter plots in $(x,\theta)$ space. Second row: time series of order parameters respectively. Third row: $(x(t),\theta(t))$ of a single swarmalator for each case. (a,d,g) Simulation parameters: $(p,K_p,K_n,T,dt,N)$=(0.75,0.1,-0.7,1000,0.25,500). Swarmalators execute swaying in four clusters.  (b,e,h) Simulation parameters: $(p,K_p,K_n,T,dt,N)$=(0.6,1,-0.5,1000,0.25,500). Swarmalators execute swirling in circular motions with noisy $S_{max}$ and $S_{min}$. (c,f,i) Simulation parameters: $(p,K_p,K_n,T,dt,N)$=(0.5,2,-0.1,1000,0.25,500). Swarmalators execute shear flow as denoted by the black arrows, in which order parameters have noisy oscillations. 
       \label{fig:unsteady}
    \end{figure*}

\textbf{1. Static Synchrony.} The swarmalators ultimately synchronize at two fixed points $(x^*,\theta^*)$ and $(x^*+\pi,\theta^*+\pi)$, where the two groups are spaced $\pi$ units equally apart. Figure~\ref{fig:Static Synchrony} shows this state where swarmalators are depicted as colored dots moving around the unit circle. The color represents the phase $\theta_i$, and the location on the circle represents the swarmalator's position $x_i$ (recall the position is a circular variable $x_i \in \mathbb{S}^1$). The rainbow order parameters $W_{\pm}$ are  plotted as larger dots (recall these are complex numbers with magnitude $< 1$ so they lie inside the unit disk) and are colored red and blue respectively, so as to distinguish them from each other (so the color does not refer to a phase, as it does for the individual swarmalators). Looking at Figure~\ref{fig:Static Synchrony}, you can see the individual swarmalators sit at fixed points with the same phase/color, and that $S_{\pm} = 1$, as expected in the static sync state. Figure~\ref{fig:xtheta_Static Synchrony} shows an alternate representation of the state: a scatter plot of the swaramlators in $(x,\theta)$ space, where conformists are colored blue, and contrarians are colored red. This sync state occurs in the limit case when all swarmalators are conformists $p=1$. This can be seen in Figure~\ref{fig:order parameters} which plots our three order parameters $S_{\text{max}}, S_{\text{min}}, V$ versus $p$ for different values of $(K_p, K_n)$. Notice in each panel the sync state is achieved when $p=1$; $S_{\pm} = 1$ and $V = 0$ at the right hand edge of the $p$-axis. Figure~\ref{fig:meshgrid} clearly illustrates the numerical and theoretical consistency in the transition between the static sync state and other states. This state was reported before \cite{o2022collective}. 
 
\textbf{2. Polarized State. } Here, the swarmalators segregate into contrarian and conformist clumps, spaced a distance of $\pi/2$ from each other as seen in Figure~\ref{fig:Polarized State}, Figure~\ref{fig:xtheta_Polarized State}.  The conformists neighbour the contrarians, and vice versa. So in a row: Conformist $\rightarrow$ Contrarian $\rightarrow$ Conformist $\rightarrow$ Contrarian. Since the conformists and the contrarians are maximally separated in their `opinions'  we call this the polarized state. The fixed points are $(x^*,\theta^*) , (x^*+\pi/2,\theta^* \pi/2), (x^*+\pi, \theta^*+\pi), (x^*+3\pi/2, \theta^*+3\pi/2)$. The order parameters take values $V=0$, $S_{max}=1$ and $0<S_{min}<1$ as illustrated in Figure~\ref{fig:order parameters} and Figure~\ref{fig:meshgrid}. This state has not been seen before in systems of swarmalators.
     
\textbf{3. Static Phase Wave.} Sometimes the swarmalators arrange themselves in a phase wave with $x_i = \pm \theta_i + C$ where the $\pm$ occur with equal probability as shown in Figure~\ref{fig:Static Phase Wave} and Figure~\ref{fig:xtheta_Static Phase Wave}. The order parameters are $V=0$, $S_{max}=1$ and $S_{min} = 0$ as seen in Figure~\ref{fig:order parameters}. Figure~\ref{fig:meshgrid} shows the region of occurrence of the state. This state was previously reported \cite{o2022collective}. 

\textbf{4. Static Asynchrony.} A static async state can be formed as well, depicted in  Figure~\ref{fig:Static Asynchrony}. It  is more clearly seen in the ($x,\theta$) plane in Figure~\ref{fig:xtheta_Static Asynchrony}.  Swarmalators are distributed uniformly, which means every phase can occur everywhere, resulting in all colors appearing everywhere, as shown in Fig.\ref{fig:Static Asynchrony}. Since $x_i$ and $\theta_i$ are uncorrelated, $S_{max}=S_{min}=0$ and $V=0$; see Figure~\ref{fig:order parameters}. {Figure~\ref{fig:meshgrid} presents the region of occurrence of the state given the numerical conditions and theoretical prediction. This state was reported before \cite{o2022collective}.

\textbf{5. Breathing polarized state.} Here the polarized state destabilizes and begins to breath as shown in Figure~\ref{fig:breathps} (Note for this and the other unsteady states we do not show the plots of them are colored dots on the unit circle; this representation was not informative. Moreover, in Figure~\ref{fig:order parameters} of $S_{\pm}(p), V(p)$ we group all these states under the umbrella ``unsteady"). The swarmalators stay in their contrarian / conformist clumps (by clump we mean a delta function mass; they all have the same position / phase) but now the clumps move in small loops about their former fixed points as indicated by the black arrows. Correspondingly, $S_{\pm}$ barely execute oscillations about their mean values (Figure~\ref{fig:xy_breathps}). However, $V>0$. Figure~\ref{fig:1_breathps} shows that the oscillations of $x(t)$ and $\theta(t)$ rise and fall incessantly and periodically. We devised two additional parameters $\delta_x$ and $\delta_\theta$, which we refer to as rotation, to distinguish among the three unsteady states, as shown in Figure~\ref{fig:rotation}. If a swarmalator completes a full rotation (from 0 to $2\pi$) in either $x$ or $\theta$, we denote $\delta_x=1$ or $\delta_\theta=1$ respectively; otherwise, they are set to 0. Figure~\ref{fig:rotation} presents the fractions of rotations of all swarmalators. Both fractions of total $\delta_x$ and $\delta_\theta$ take on the value of 0.  The state is novel.

 \textbf{6. Swirling} In this state the conformists stay in their clumps, but the contrarian break out into a noisy vortex like structure as seen in Figure~\ref{fig:swirling}; note the red dots are dispersed, but the two blue clumps remain. The vortices periodically form and disperse, and within each vortex the contrarians swirl as indicated by the black arrows. This vacillatory motion manifests as irregular times series of $S_{\pm}$ as shown in Figure~\ref{fig:xy_swirling}. Fractions of $\delta_x=1$, while fractions of $\delta_\theta<1$. By combining the oscillatory behavior in Figure~ \ref{fig:1_swirling} with the fractions of rotations in Figure~\ref{fig:rotation}, the state can be distinguished. This state has not been reported before. 
  
\textbf{7. Active bands} Swarmalators form band like structure in $(x,\theta)$ space which move and break up periodically (Fig~\ref{fig:active_bands}). Here the swarmalators' positions are almost synchronized, but their phases are distributed. The overall macroscopic motion is somewhat irregular as indicated by the time series of $S_{\pm}$ (Fig~\ref{fig:xy_active_band}). Figure~\ref{fig:1_active_band} shows the oscillatory behavior of a typical swarmalator (the vortices appear and disperse periodically). Both fractions of $\delta_x$ and $\delta_\theta$ are equal to 1. By combining the pattern in Figure~\ref{fig:1_active_band} with rotations, we can determine this state, as shown in Figure~\ref{fig:rotation}. This state is to our knowledge novel.

         \begin{figure}
    \centering
    \includegraphics[width= 1 \columnwidth]{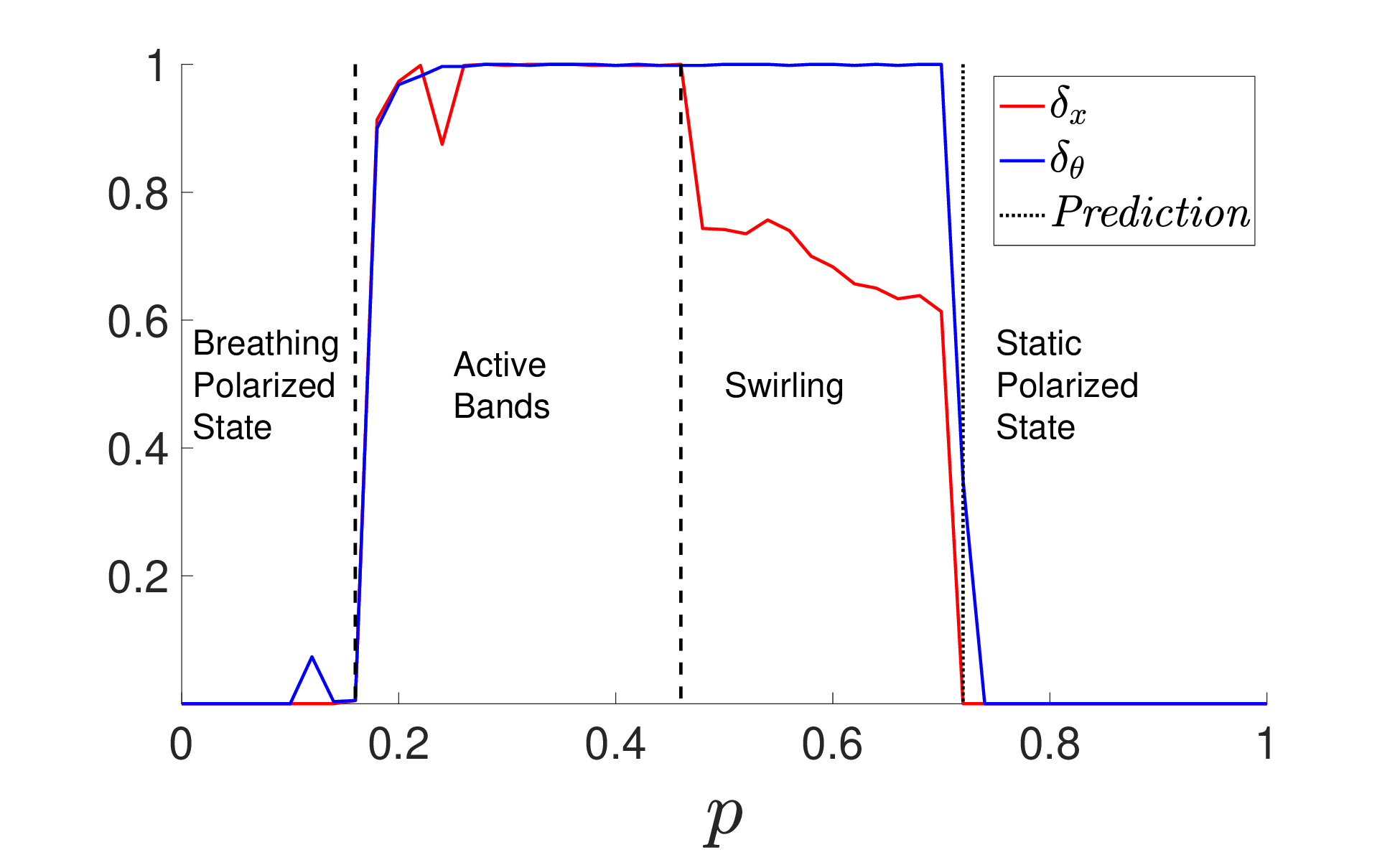}  
    \caption[ ]
    {\small \textbf{Fractions of rotations  plot for different $p$.}} $\delta_x$, $\delta_\theta$ represent the rotation of $x$ and $\theta$ respectively. The dotted line is the theoretical prediction of the analysis of static polarized state. ($K_p,K_n$) = (1.5,-0.6). $(T,N,J,dt)$=$(2000,600,1,0.1)$. First $30\%$ data are dropped.  
        \label{fig:rotation}
\end{figure}

\section{Analysis}

Here we analyze the stability of the static async and polarized states. The static phase wave, though static, was too difficult to crack. The same is true of the breathing polarized, swirling, and active band states; being unsteady, their analysis was intractable.


\textbf{Polarized State}. We analyze the stability of the state for all finite $N$ using standard methods, namely linearizing around the fixed points: $(x_i,\theta_i) = (x^*,\theta^*), (x^*+\pi/2, \theta^*+\pi/2), (x^*+\pi, \theta^*+\pi), (x^*+3\pi/2,\theta^*+ 3\pi/2)$. The algebra is somewhat involved, but the essence of our approach is simple: take advantage of the block structure of the associated Jacobean $M$. It turns out its easier to move to $(\xi, \eta)$ coordinates defined by $\xi_i = x_i + \theta_i, \eta_i = x_i - \theta_i$. The ODEs in this frame are
\begin{align}
    \dot{\xi_i} &=  J_+ S_+ \sin( \Phi_+ - \xi ) + J_{-} S_{-} \sin( \Phi_- - \eta ) \label{eom_xi}, \\
    \dot{\eta_i} &=  J_+ S_+ \sin( \Phi_+ - \xi ) + J_+ S_- \sin( \Phi_- - \eta ), \label{eom_eta} 
\end{align}
where $J_\pm =(J\pm K)/2.$

\noindent
The Jacobian has form
\begin{equation}
    M = \left[ 
\begin{array}{cc} 
  Z_{\xi} & Z_{\eta} \\ 
  N_{\xi} & N_{\eta} \\
\end{array} 
\right] \label{eq:Jacobian}
\end{equation}
where
\begin{align}
(Z_{\xi})_{ij} &= \frac{\partial \dot{\xi}_i}{\partial \xi_j} \\
(Z_{\eta})_{ij} &= \frac{\partial \dot{\xi}_i}{\partial \eta_j} \\
(N_{\xi})_{ij} &= \frac{\partial \dot{\eta}_i}{\partial \xi_j} \\
(N_{\eta})_{ij} &= \frac{\partial \dot{\eta}_i}{\partial \eta_j} 
\end{align}
Plugging the values for the derivatives yields at the fixed points gives
\begin{equation}
    M_{ps} = \left[ 
\begin{array}{cc} 
  A(K_p, K_n) & B(-K_p, -K_n) \\ 
  A(-K_p, K_n) & B(K_p, K_n) \\
\end{array} 
\right] 
\end{equation}
where $A$ is
\[   
A_{ij}(x,y) = 
     \begin{cases}
       \frac{-(N-1)}{2N}(J + x) & i = j, i < n_{p}\\
       \frac{-(N-1)}{2N}(J + y) & i = j, i \geq  n_{p} \\
       \frac{J + x}{2 N} & i \neq j, i < n_{p}\\
       \frac{J + y}{2 N} & i \neq j, i \geq  n_{p} \label{eq:Aij}
     \end{cases}
\]
and $B$ is
\[   
B_{ij}(x,y) = 
     \begin{cases}
       \frac{2(N/2 - n_p-1) + 1)}{2N} (J + x) & i = j, i < n_{p}\\
       \frac{2(N/2 - n_p-1) + 1)}{2N} (J + y)  & i = j, i \geq  n_{p} \\
       \frac{J + x}{2 N} & i < n_{p}, i < n_p\\
       -\frac{J + x}{2 N} & i < n_{p}, i \geq n_p\\
       -\frac{J + y}{2 N} & i \geq n_{p}, i < n_p\\
       \frac{J + y}{2 N} & i \geq n_{p}, i \geq n_p\\ \label{eq:Bij}
     \end{cases}
\]
where $ n_{p} = ceil(p N)$ is the number of swarmalators with $K = K_p$. Intuitively, what is going on here is that each $A,B$ are subdivided into contrarian and conformist populations. Look at the diagonal elements of $A$: the first $n_p$ have coupling $K_p$, while the remaining $N-n(p)$ have coupling $K_n$. We write $A,B$ below for the $(n,p) = (4,1/4)$ so the structure can be seen visually.

\[
A =\left[
\begin{array}{cccc}
 -\frac{3}{8} (x+1) & \frac{x+1}{8} & \frac{x+1}{8} & \frac{x+1}{8} \\
 \frac{y+1}{8} & -\frac{3}{8} (y+1) & \frac{y+1}{8} & \frac{y+1}{8} \\
 \frac{y+1}{8} & \frac{y+1}{8} & -\frac{3}{8} (y+1) & \frac{y+1}{8} \\
 \frac{y+1}{8} & \frac{y+1}{8} & \frac{y+1}{8} & -\frac{3}{8} (y+1) \\
\end{array}
\right]
\]

\[
B =\left[
\begin{array}{cccc}
 \frac{3 (x+1)}{8} & \frac{1}{8} (-x-1) & \frac{1}{8} (-x-1) & \frac{1}{8} (-x-1) \\
 \frac{1}{8} (-y-1) & \frac{1}{8} (-y-1) & \frac{y+1}{8} & \frac{y+1}{8} \\
 \frac{1}{8} (-y-1) & \frac{y+1}{8} & \frac{1}{8} (-y-1) & \frac{y+1}{8} \\
 \frac{1}{8} (-y-1) & \frac{y+1}{8} & \frac{y+1}{8} & \frac{1}{8} (-y-1) \\
\end{array}
\right]
\]

Now, getting back to our goal of finding the eigenvalues $\lambda$ of $M$. A well known fact for block matrices is $\det(M) = \det(AD - BC)$ if the sub-matrices $AD, BC$ commute, which is the case for us. To find the $\lambda$, we write $\det(M - \lambda I) = \det( (A-\lambda I) (D - \lambda I)) - BC)$. We define 
\begin{align}
    G := (A-\lambda I) (D - \lambda I)) - BC
\end{align}
G inherits the structure of $B$ and has thus $6$ unique elements:
\[   
G = 
     \begin{cases}
       g_1(p,N,\lambda)  & i = j, i < n_{p}\\
       g_2(p,N,\lambda)  & i = j, i \geq  n_{p} \\
       g_3(p,N,\lambda) & i < n_{p}, i < n_p\\
       g_4(p,N,\lambda) & i < n_{p}, i \geq n_p\\
       g_5(p,N,\lambda) & i \geq n_{p}, i < n_p\\
       g_6(p,N,\lambda) & i \geq n_{p}, i \geq n_p
     \end{cases}
\]
Or in block format:
\begin{equation}
    G = \left[ 
\begin{array}{cc} 
  G_1 & G_2 \\ 
  G_3 & G_4 \\
\end{array} 
\right] 
\end{equation}
where $(G_{1})_{i,j} = g_1, \dots $. We want to find the eigenvalues $\hat{\lambda}$ of $G$ for which we need $\det (G - \hat{\lambda} I)$. Notice, however, that the sub-matrices $G_i$ are non-square, so we can't use the previous formula we used $\det(M) = \det(AD-BC)$. Instead we use Schur's formula:
\begin{align}
\det (G - \hat{\lambda} I) & = \det(G_1 - \hat{\lambda} I ) \det ((G_4 - \hat{\lambda}I) \nonumber \\
& - G_3 ( G_1 - \hat{\lambda})^{-1} G_2) \\
\det (G - \hat{\lambda} I) & =  \det(G_1 - \hat{\lambda} I ) \det (G_5)
\end{align}
Now all that's left is to find expressions for the determinants. This was the bottle neck in the calculation. After much algebra, we find

\begin{align}
    \det(G_1 - \hat{\lambda} I ) &= \left(g_1-g_3\right){}^{n_p-1} \left(g_3 \left(n_p-1\right)+g_1\right) \nonumber \\
    & -\hat{\lambda } \left(\left(g_1-g_2\right) n_p^2\right){}^{n_p-1} \\
    \det(G_5) &= \frac{\left( \tilde{a} +\tilde{b} \hat{\lambda} + \hat{\lambda }^2\right) \left(g_6-g_2+\hat{\lambda }\right){}^{n_q-1}}{g_3 \left(n_p-1\right)+g_1-\hat{\lambda }}
\end{align}
where $n_q := N - n_p$ is the number of swarmalators with $K = K_n$ and $\tilde{a} = \sum_j \tilde{a}_i g_i$ and $\tilde{b} = \sum_{i,j} \tilde{b_{i,j}} g_i g_j$ (for convenience we do not write out $\tilde{a_i}, \tilde{b_i}$ ). Multiplying these together and equating to zero yields four distinct eigenvalues:
\begin{align}
\hat{\lambda}_0(p,N,\lambda) &=  g_1 - g_3 \hspace{0.25 cm} w.m \; \; n_p - 1 \\
\hat{\lambda}_1(p,N,\lambda)  &= g_2 - g_6  \hspace{0.25 cm} w.m \; \; n_q - 1 \\
\hat{\lambda}_{2/3}(p,N,\lambda)  & =  \frac{1}{2} \sum_i a_i g_i \\ 
 & \pm \frac{1}{2} \sqrt{ \Big( \sum_i b_i g_i \Big)^2 - 4  \sum_{i,j} c_{ij} g_i g_j  } \hspace{0.25 cm} w.m \; \; 1 
\end{align}
where \( w.m \) means 'with multiplicity', and we have dropped the dependence on \(\lambda\) for the \(g_i\). The other coefficients depend on \(p,N\): \(a_i = a_i(p,N)\), \(b_i = b_i(p,N)\), and \(c_{i,j} = c_{i,j}(p,N)\). Now, recall these \(\hat{\lambda}\) are the eigenvalues of \(G\), but our target are those of \(M\). So we set \(\hat{\lambda}_i = 0\) and solve for \(\lambda\). After much calculation, we eventually derive:
\begin{align}
    & \lambda_0 = 0 \\
    & \lambda_1 = \pm \frac{1}{2} \sqrt{ J(q K_p + p K_n)}  \\
    & \lambda_2 = \frac{1}{2} \Big( ( -p(J + K_p) \pm \sqrt{p^2(J+K_p)^2 - 8 J K_p p + 4 J K_p } \Big)  \\
    & \lambda_3 = \frac{1}{2} \Big( -q(J + K_n) \pm \sqrt{q^2(J+K_n)^2 - 8 J K_n q + 4 J K_n } \Big)
\end{align}
where $q = 1-p$ and the multiplicities are $2, 1, n_p-1, n_q-1$. Notice that $\lambda_{2}$ becomes $\lambda_3$ under the transformation $(p,K_p) \rightarrow (q, K_n)$. 

We've done the hard work. Now its time to use the $\lambda$ to deduce the stability of the polarized state. The zeroth $\lambda_0$ does not play a role in the bifurcation; it simple corresponds to the rotational symmetry in the model. The first $\lambda_1$ however undergoes a (degenerate) saddle node bifurcation when the argument of the radical becomes real. The second $\lambda_2$ has negative real parts for all parameter regimes of interest and so is unimportant, while the final $\lambda_3$ undergoes a (degenerate) hopf bifurcation. Thus, the state is stable when
\begin{align}
    J(q K_p + p K_n) < 0, \\
    J + K_n > 0, \\
    p \ge \frac{1}{2}.
\end{align}
When $J=1$, the critical $K_n$ values are:
\begin{align}
    K_{n} &= -\frac{q}{p} K_p, \\
    K_{n} &= -1.
\end{align}
The dotted black lines in Figure~\ref{fig:order parameters} shows these predictions are consistent with numerics.

\textbf{Static Sync State}. The fixed points of this state are $(x^*,\theta^*)$ and $(x^*+\pi, \theta^*+\pi)$. Similarly, by linearizing around the fixed points in ($\xi,\eta$) space. We seek the eigenvalues $\lambda$ of the Jacobian $M$ in Eq.~\ref{eq:Jacobian}. Plugging the values for the derivatives yields at the fixed points gives \begin{equation}
    M_{ss} = \left[ 
\begin{array}{cc} 
  A(K_p, K_n) & A(-K_p, K_n) \\ 
  A(-K_p, K_n) & A(K_p, K_n) \\
\end{array} 
\right] ,
\end{equation}where $A$ is the same subblock in Eq.~\ref{eq:Jacobian}.}

By using the following identity for symmetric block matrices, the egienvalues $\lambda$ of $M_{ss}$ can be found:
\begin{equation}
    \det E := \left[ 
\begin{array}{cc} 
  C & D \\ 
  D & C\\
\end{array} 
\right] = \det (C+D) \det (C-D).
\end{equation}
Applying this identity to $M_{ss}$ yields:
\begin{align}
    & \lambda_0 = 0, \\
    & \lambda_1 = - J,  \\
    & \lambda_2 = -\frac{K_p (N-n_p)+K_nn_p}{N},  \\
    & \lambda_3 = -K_p, \\
    & \lambda_4 = -K_n,
\end{align} 
with  multiplicities $2, N-1, 1, floor(n_p/2), N-2-floor(n_p/2)$. Hence, the state is stable when\begin{align}
    J&> 0, \\
    K_n,K_p&>0,
\end{align}
which also means $p=1$. The dashed black lines in Figure~\ref{fig:order parameters} indicate that these predictions are consistent with the numerical results.

\textbf{Incoherence}. This state is analyzed is the same manner as previous studies: we take the $N \rightarrow \infty$ limit and perturb around $\rho_0 = (2 \pi)^{-2}$. The calculation is virtually the same, and the result is the same also, so we just quote the result:
\begin{align}
    & \langle K \rangle_c = - \langle J \rangle  
\end{align}
where the $\langle . \rangle$ denotes the average. For the $h(K) = p \delta(K-K_n) + q \delta(K+K_p)$ example this becomes
\begin{equation}
p_s = \frac{1 + K_n}{K_n - K_p}
\label{p_s}
\end{equation}
which, interestingly, is identical to the result found for constant coupling \cite{o2022collective} and $K_j$ coupling \cite{o2022swarmalators}. Figure~\ref{fig:order parameters} shows these predictions are consistent with numerics.

\section{Other Coupling Distributions}
\begin{figure}
        \centering
        \begin{subfigure}[b]{0.23\textwidth}
           \centering
            \includegraphics[width=\textwidth]{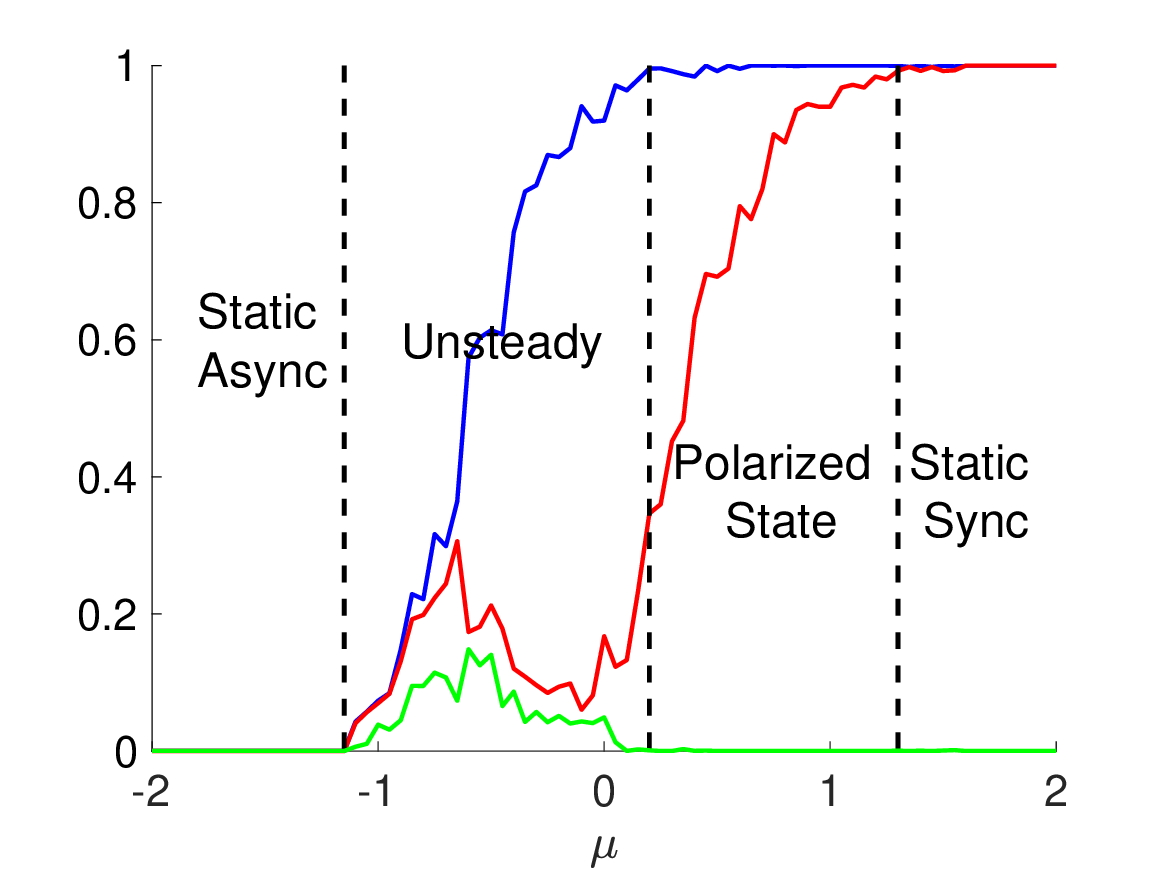}
            \caption[Single Gaussian]%
            {{\small Single Gaussian}}    
            \label{fig:single_gaussian}
        \end{subfigure}
        \hfill
        \begin{subfigure}[b]{0.23\textwidth}  
            \centering
            \includegraphics[width=\textwidth]{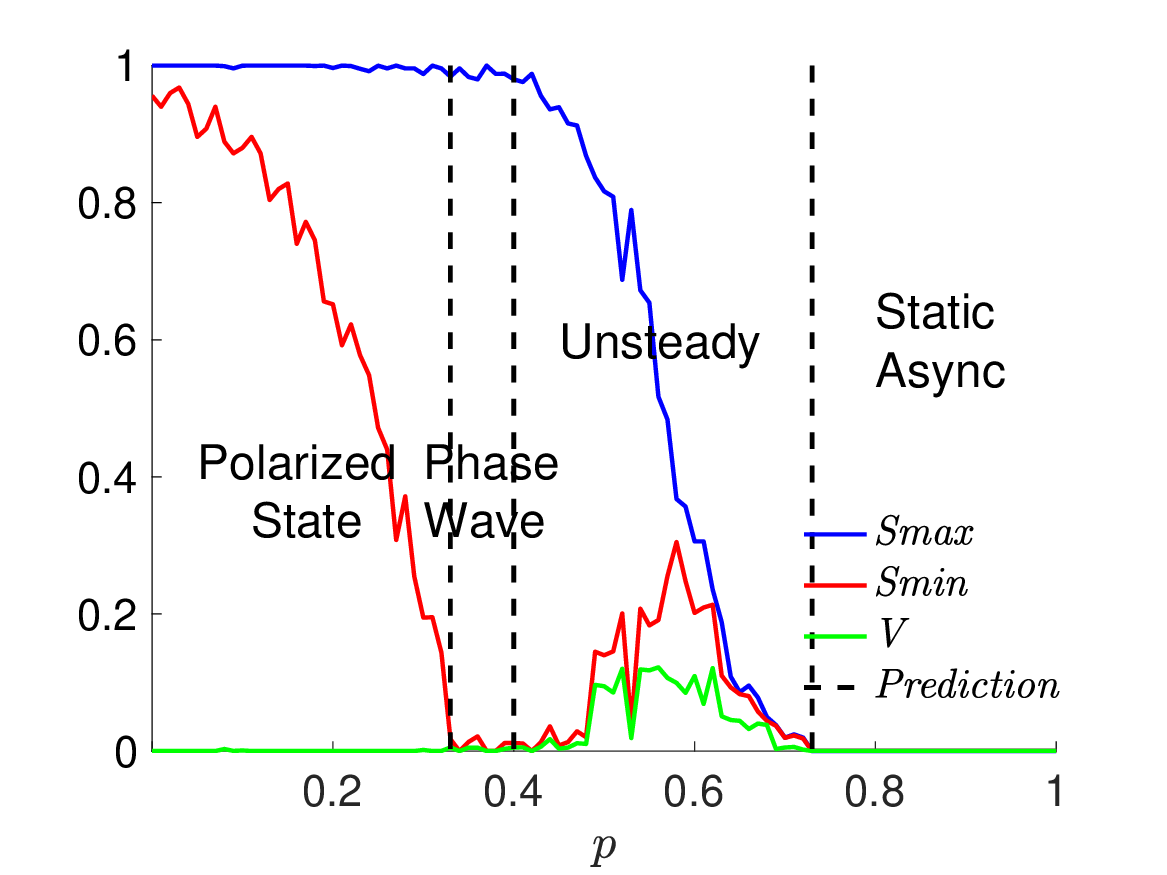}
            \caption[Mixed Gaussian]%
            {{\small Mixed Gaussian}}    
            \label{fig:mixed_gaussian}
        \end{subfigure}
        \caption[ other distributions ]
        {\small \textbf{Order parameters and averaged velocity plots
for other coupling distributions.} The same structure is observed (a) Gaussian with variance $\sigma=0.5$. (b) Mixture of gaussians with $(K_p,K_n,\sigma) = (-2, 1, 0.5).$ Simulation parameters: $(J,T,dt,N)=(1,2000,0.5,500).$} 
        \label{fig:other distributions}
    \end{figure}
Besides the double delta distribution with certain fractions $p$ and $q$, we also studied the phase couplings $K$ from other distributions: \begin{itemize}
    \item[1] Single Gaussian Distribution: $h(K) \sim N(\mu,\sigma|K)$,
    \item[2] Mixed Gaussian Distribution: $h(K) \sim pN(K_n,\sigma|K)+(1-p)N(K_p,\sigma|K)$,
\end{itemize}where $N(\mu,\sigma|x)$ is the normal distribution for random variable $x$. The same states were also found in each case. Figure~\ref{fig:other distributions} summarizes them by plotting the order parameters. Besides, By setting $(\mu,\sigma)=(-0.5,0.1)$, the system with single Gaussian coupling distribution gives us a static phase wave state which is not shown in Figure~\ref{fig:single_gaussian}. Also, we can get static sync by letting $(p,K_p,\sigma)=(0,2,0.1)$ in the mixed Gaussian case.

\section{Match to Real-world Swarmalators}
Sperm are prototypical microswimmers that exhibit collective behavior by synchronizing their tail movements while swarming in a solution \cite{yang2008cooperation}. When contained within 1D rings, sperm obtained from ram semen transition from an isotropic state resembling the static async state to a vortex state, where sperm rotate either clockwise or counterclockwise \cite{creppy2016symmetry}. This suggests that their positions and orientations are arranged in a manner akin to the static phase wave. It's worth noting that the static phase wave actually represents a state of uniform rotation and remains motionless in the frame that moves along with the natural frequencies $\omega,\nu$ \cite{creppy2016symmetry}. Additionally, the transition involves a temporary reduction in rotational velocity, as observed in the decay shown in Figure 4(a) of reference \cite{creppy2016symmetry}, which aligns with characteristics of a Hopf bifurcation, just similar to the ring model. It is important to highlight that, unlike other research on synchronizing sperm, in this case, the phase variable refers to the orientation of the sperm \cite{yang2008cooperation}, not their tail rhythm. 

Vinegar eels, a type of nematode found in beer mats and tree wound slime \cite{quillen2021metachronal,peshkov2022synchronized}, are swarmalators because they sync the wriggling of their heads and swarm in solution. Their collective movement suggests a potential interaction between this synchronization and the swarming behavior \cite{quillen2021metachronal,peshkov2022synchronized} (neighboring eels synchronize more easily than those at a distance, indicating an interaction between synchronization and swarming. The synchronized eels likely influence their local hydrodynamic environment, subsequently affecting each other's movements, thus showcasing an interaction between swarming and synchronization). When these eels are confined to 2D disks and maneuver near the 1D ring boundary, they create metachronal waves. In these waves, the pattern of their gait phase and their spatial positions around the ring resemble the configuration seen in the static phase wave \cite{quillen2021metachronal,peshkov2022synchronized}. However, it's essential to note that the metachronal waves possess a winding number $k > 1$, signifying that a full rotation in physical space $x$ leads to $k > 1$ rotations in phase $\theta$.

\section{Discussion}
Swarmalators are a new sub-field with little to no theoretical results. Our work is as part of a research series \cite{o2022collective,yoon2022sync,o2022swarmalators} whose goal is to develop a theory for swarmalators by first focusing on the simplest models possible (the 1D model presented here) and trying to solve those. That is, we follow the the minimal modeling or physicists tradition. This is the approach Winfree and later Kuramoto took with their famous coupled oscillator models which essentially launched the field; we follow in their footsteps.

The naked 1D swarmalator model, that with uniform natural frequencies $(\nu,\omega)$ and couplings $(J,K)$, was first introduced and solved in \cite{o2022collective} Then distributed natural frequencies $(\nu, \omega) \rightarrow (\nu_i, \omega_i)$ \cite{yoon2022sync} and random $K_j$-couplings $(J,K) \rightarrow (J_j, K_j)$ \cite{o2022swarmalators} were studied. Here we tackled $K_i$ couplings $(J,K) \rightarrow (J_i, K_i)$ and found some new states: the polarized state, which we characterized analytically, as well as the unsteady breathing, swirling, and active band states, which we characterized numerically. In the future work, before delving into the study of van Hemmen couplings $(J,K) \rightarrow (J_{ij}, K_{ij})$, it would also be interesting to study the role that $J$ plays. We will remove the simplicity $J=1$ and investigate the consequences when this position coupling strengths are also chosen in a similar manner like $K$’s ($K_{i}$'s and $K_{ij}$'s).  Next, we will complete the random coupling sequence by studying van Hemmen couplings $(J,K) \rightarrow (J_{ij}, K_{ij})$, where $K_{ij}$ is drawn from $h(K) = \frac{1}{4}\delta(K-\mu+\gamma)+\frac{1}{2}\delta(K-\mu)+\frac{1}{4}\delta(K-\mu-\gamma)$ (and in theory the same for $J$; although we may set $J_{ij} = 1$ for simplicity) This is a crude two parameter representation of a unimodal distribution with mean $\mu$ and $\gamma$. Van Hemmen couplings $K_{ij}$ are well studied in the Ising model of statistical physics, and are more realistic than the $K_i, K_j$ couplings models (since then randomness is associated with an interaction between a pair $(i,j)$ of oscillators, which is more common in nature) and thus may be applicable to real world swarmalators such as vinegar eels and sperm. 


\section{Acknowledgement}
This research was partially supported by NSF-$2225507$ (M.Z.).

\bibliographystyle{apsrev}

\end{document}